\documentclass{svproc}

\usepackage{amsmath,amssymb,graphics,epsfig,color,cite}
\usepackage{bm}
\newcommand {\pa}{\partial}
\usepackage{tikz}
\def\Ran{{\rm Ran}}
\usetikzlibrary{arrows}

\usepackage{lmodern}
\usepackage{amsmath, amssymb}
\newcommand{\indep}{\rotatebox[origin=c]{90}{$\models$}}
\usepackage{url}

\begin{document}
\mainmatter              
\title{Exit event from a metastable state and Eyring-Kramers law  for the overdamped Langevin dynamics}
\titlerunning{Metastability and Eyring-Kramers law}  
%

\author{Tony Leli\`evre\inst{1} \and Dorian Le Peutrec\inst{2}
 \and Boris Nectoux\inst{1}}
\authorrunning{Tony Leli\`evre et al.} 
%
%
\institute{\'Ecole des Ponts, Universit\'e  Paris-Est, INRIA, 77455 Champs-sur-Marne, France,\\
   \email{ \{tony.lelievre,boris.nectoux\}@enpc.fr} 
\and
Laboratoire de Math\'ematiques d'Orsay, Univ. Paris-Sud, CNRS, Universit\'e Paris-Saclay, 91405 Orsay,  \email{ dorian.lepeutrec@math.u-psud.fr} }
%

\maketitle              

\begin{abstract}
\begin{sloppypar}
\noindent
In molecular dynamics,  several  algorithms have been designed over the  past few years to  accelerate the sampling of  the exit event from a metastable domain, that is to say the time spent and the exit point from the domain. Some of them  are based on the fact that the exit event from a metastable region is well approximated by a Markov jump process. In this work, we present recent results on the exit event from a metastable region  for the overdamped Langevin dynamics obtained in~\cite{BN2017,di-gesu-lelievre-le-peutrec-nectoux-16,DLLN}. These results aim in particular at justifying  the use of  a Markov jump process parametrized by  the Eyring-Kramers law to model the  exit event from a metastable region.

\end{sloppypar}
\end{abstract}
The objective of this note is to give motivations (Section~\ref{section:1}) and outlines of the proofs (Section~\ref{sec.mainresults}) of results recently
obtained in\cite{BN2017,di-gesu-lelievre-le-peutrec-nectoux-16,DLLN}. These results justify the use of the Eyring-Kramers formulas together with a kinetic Monte Carlo model
to model the exit event from a metastable state for the overdamped Langevin dynamics. Such results are particularly
useful to justify algorithms and models which use such formulas to build reduced description of the overdamped
Langevin dynamics.
\section{Exit event from a metastable domain and  Markov jump process} 
\label{section:1}
\subsection{Overdamped Langevin dynamics and metastability} 
 \label{sec11}
Let~$(X_t)_{t\ge0}$ be the stochastic process  solution to the  overdamped Langevin dynamics in~$\mathbb R^d$:
\begin{equation}\label{eq.langevin}
d X_t = -\nabla f(X_t) d t + \sqrt{h} \ d B_t,
\end{equation}
where~$f\in C^{\infty}(\mathbb R^d,\mathbb R)$ is the potential function,~$h>0$ is the temperature and~$(B_t)_{t\geq 0}$ is a standard~$d$-dimensional Brownian motion. The overdamped Langevin dynamics can be used  for instance to describe the motion of the atoms of a molecule or the diffusion of impurities in a crystal (see for instance~\cite[Sections 2 and 3]{MaSc} or \cite{chandrasekhar1943stochastic}). The term $-\nabla f(X_t) $ in~\eqref{eq.langevin} sends the process towards local minima of~$f$, while thanks to the noise term $\sqrt{h} \ d B_t$, the process  $X_t$ may jump from one basin of attraction of the dynamics $\dot{x}=-\nabla f(x)$ to another one. If the temperature is small (i.e. $h \ll1$), the process $(X_t)_{t\ge 0}$ remains during a very long period of time trapped around a neighborhood of a local minimum of~$f$, called a metastable state, before going to another region. For that reason, the process~\eqref{eq.langevin} is said to be metastable.   
More precisely, a  domain~$\Omega\subset \mathbb R^d$ is said to be  metastable for the probability measure~$\mu$ supported in $\Omega$ if, when $X_0\sim \mu$, the process \eqref{eq.langevin} reaches a local equilibrium  in $\Omega$ long before escaping from it. This will be made
more precise below
using the notion of
quasi-stationary
distribution (see Section~\ref{sec.qsd}). The  move from one metastable region  to another is typically related to a macroscopic change of  configuration of the system. Metastability implies  a  separation of timescales which is one of the major issues when trying to have access to the macroscopic evolution of the system  using simulations made at the microscopic level.
 Indeed, in practice,    many transitions cannot be observed by integrating directly the trajectories of the process~\eqref{eq.langevin}.
 To overcome this difficulty,  some algorithms use   the fact that the exit event from a metastable region can be well approximated by a Markov jump process with transition rates computed with the Eyring-Kramers formula, see for example the Temperature
Accelerated Dynamics method~\cite{sorensen-voter-00} that will be described below.

 \subsection{Markov jump process and Eyring-Kramers law} 
\label{sec.JP}
\subsubsection{Kinetic Monte Carlo methods.} 
Let $\Omega\subset \mathbb R^d$ be a domain of the configuration space and let us assume that the process~\eqref{eq.langevin} is initially distributed according to the probability measure $\mu$ (i.e. $X_0\sim\mu$) which is supported in $\Omega$ and for which the exit event from~$\Omega$ is metastable. Let us denote by~$(\Omega_i)_{i=1,\ldots,n}$ the  surrounding domains of~$\Omega$ (see Figure~\ref{fig:surrounding}), each of them corresponding to a macroscopic state of the system.
Many reduced
models and
algorithms rely on
the fact that the exit event from $\Omega$, i.e. the next visited state by the process~\eqref{eq.langevin} among the~$\Omega_i$'s as well as the time spent by the process~\eqref{eq.langevin} in $\Omega$, is efficiently approximated by a Markov jump process using kinetic Monte Carlo methods~\cite{schuette-98, schuette-sarich-13, voter-05,wales-03,cameron-14b,fan-yip-yildiz-14}. Kinetic Monte Carlo methods  simulate a Markov jump process in a discrete state space.  
To use a kinetic Monte Carlo algorithm in order to sample the exit event from~$\Omega$, one needs   for~$i\in\{1,\ldots,n\}$  the transition rate $k_i$  to go from the state~$\Omega$ to the state~$\Omega_i$. A kinetic Monte Carlo algorithm  generates the next visited state $Y$ among the~$\Omega_i$'s and the time $T$ spent  in $\Omega$ for  the process~\eqref{eq.langevin}  as follows: 
\begin{enumerate}
\item First sample   $T$  as an exponential random variable with parameter   $\sum_{i=1}^n k_i$, i.e.:
\begin{equation}\label{eq.Te}
T\sim \mathcal E\Big ( \sum_{i=1}^n k_i\Big).
\end{equation}
\item Then, sample the next visited state~$Y$ independently from $T$, i.e
\begin{equation}\label{eq.Ye}
Y \  \indep \ T
\end{equation}
 using the following law : for all $i\in\{1,\ldots,n\}$, 
 \begin{equation}\label{eq.Ylaw}
 \mathbb P\big[  Y=\,  i  \big]\ =\ \frac{k_{i}}{  \, \sum_{\ell=1}^n k_{\ell}}.\end{equation}
\end{enumerate} 
 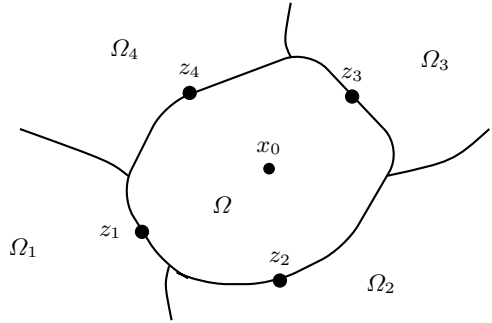
\begin{figure}
\begin{center}
\begin{tikzpicture}[scale=0.48]
\tikzstyle{vertex1}=[draw,circle,fill=black,minimum size=4pt,inner sep=0pt]
\tikzstyle{vertex}=[draw,circle,fill=black,minimum size=5pt,inner sep=0pt]
\tikzstyle{ball}=[circle, dashed, minimum size=1cm, draw]
\tikzstyle{point}=[circle, fill, minimum size=.01cm, draw]
\draw [thick, rounded corners=10pt] (1,0.5) -- (-0.25,2.5) -- (1,5) -- (5,6.5) -- (7.6,3.75) -- (6,1) -- (4,0) -- (2,0) --cycle;
 \draw  (2.6,2.2) node[]{$\Omega$};
     \draw[thick]   (-3,4.3)  ..controls   (-0.5,3.4)..   (0,3) ;
      \draw[thick]   (10,4.3)  ..controls   (9,3.4)..   (7.16,3) ;
        \draw[thick]   (1.13,0.5)  ..controls   (1,0)..  (1.2,-1) ;
        \draw[thick]   (4.5,6.3)  ..controls   (4.3,6.7)..  (4.5,7.8);
   
   \draw  (3.9,3.2) node[vertex1,label={$x_0$}](v){};
\draw (1.7,5.3) node[vertex,label={$z_{4}$}](v){};
 \draw   (-0.1,6.5) node[]{$\Omega_4$};
\draw (4.2,0.1) node[vertex,label= {$z_{2}$}](v){};
 \draw   (7,0) node[]{$\Omega_2$};
\draw (0.38,1.45) node[vertex,label={}](v){};
 \draw  (-0.5,1.45) node[]{$z_{1}$};
  \draw   (-2.9,1) node[]{$\Omega_1$};
\draw (6.2,5.2) node[vertex,label= {$z_{3}$}](v){};
 \draw   (8.5,6.2) node[]{$\Omega_3$};
\end{tikzpicture}
\caption{Representation of the domain $\Omega$, the surrounding domains~$(\Omega_i)_{i=1,\ldots,4}$ of~$\Omega$, the global minimum~$x_0$ of $f$ in~$\Omega$ and $\{z_i\}={\rm argmin}_{\pa \Omega\cap \Omega_i}f$  ($i\in\{1,2,3,4\}$).  }
 \label{fig:surrounding}
 \end{center}
\end{figure}
\begin{remark}\label{eq.equi-kmc}
Let us give an equivalent way to sample $T$ and $Y$ in a Monte Carlo method.  Let $(\tau_i)_{i\in \{1,\ldots,n\}}$  be $n$  independent random variables such that for all $i\in\{1,\ldots,n\}$, $\tau_{i}$ is exponentially distributed with parameter $k_{i}$. Then, the couple $(T,Y)$ has the same law as $(\min_{ j\in \{1,\ldots,n\} }  \tau_{j},{\rm argmin}_{j\in \{1,\ldots,n\}}  \tau_{j})$. 
 \end{remark}
 \subsubsection{Eyring-Kramers law.} 
In practice, the  transition rates~$(k_i)_{i\in\{1,\ldots,n\}}$ are computed using the Eyring-Kramers formula~\cite{hanggi-talkner-barkovec-90,voter-05}: 
\begin{equation}\label{eq.ek}
 k_{i}=A_{i}\, e^{-\frac 2h   (f(z_{i})-f(x_0) )}, 
\end{equation}
where $x_0\in \Omega$ is the unique global minimum of $f$ in $\overline \Omega$ and $\{z_i\}={\rm argmin}_{\pa \Omega\cap \pa \Omega_i}f$, see Figure~\ref{fig:surrounding}. We here assume for simplicity that the minimum is attained at one single point $z_i$ but the results below can be generalized to more general settings. If $\Omega$ is the basin of attraction of $x_0$ for the dynamics $\dot x=-\nabla f(x)$  so that $z_i$ is a  saddle point of $f$ (i.e. a critical point of index $1$),  then, for the overdamped Langevin dynamics~\eqref{eq.langevin}, the prefactor $A_i$ writes:
\begin{equation}\label{eq.prefactor}
 A_{i}=\frac{\vert \lambda(z_{i})\vert}{2\pi }\frac{  \sqrt{ {\rm det\,  Hess}\, f (x_0) }    }{   \sqrt{\vert {\rm det \,Hess}\, f(z_{i})  \vert  }   },
 \end{equation}
 where $\lambda(z_{i}) $ is the negative eigenvalue of the Hessian matrix of $f$ at $z_{i}$. Notice that the formula~\eqref{eq.prefactor}  requires that $x_0$ and $z_i$ are non degenerate critical points of~$f$. 
 The formulas~\eqref{eq.ek} and~\eqref{eq.prefactor} have been first obtained in the small temperature regime by Kramers~\cite{kramers-40}  (see the review of the literature~\cite{hanggi-talkner-barkovec-90}). 
 
\begin{remark}
In the Physics literature,  the approximation of the macroscopic evolution of the system with a  Markov jump process with transition rates   computed with the Eyring-Kramers formula~\eqref{eq.ek}-\eqref{eq.prefactor} is sometimes called the Harmonic  Transition State Theory~\cite{marcelin-15,vineyard-57}. 
\end{remark}

 \subsection{The {\em temperature accelerated dynamics} algorithm.} 
 \begin{sloppypar}
 The temperature accelerated dynamics (TAD) algorithm  proposed by  M.R. S$\o$rensen and A.F~Voter~\cite{sorensen-voter-00} aims at   efficiently approximating the   exit event from a metastable domain for the dynamics~\eqref{eq.langevin}   in order to have access to the macroscopic evolution of the system. We also refer to~\cite{aristoff-lelievre-14}  for a  mathematical analysis of this algorithm in a one-dimensional setting. \end{sloppypar}
 \medskip
 
 \noindent 
 The basic idea of the   TAD  algorithm is the following: the  exit time from the metastable domain~$\Omega$ increases exponentially  with the inverse of the temperature, see indeed~\eqref{eq.Te}-\eqref{eq.ek}. The idea is then to simulate the process at higher temperature to accelerate the simulation of the exit event. 
 Let us assume that the process~$(X_t)_{t\ge 0}$, evolving at the temperature $h_{low}$ is at some time $t_0\ge 0$ in  the domain~$\Omega\subset\mathbb R^d$ which is metastable for the initial condition~$X_{t_0}\in \Omega$.
Following~\cite{sorensen-voter-00}, let us assume that the process instantaneously reaches the local equilibrium in~$\Omega$, i.e. that $X_{t_0}$ is distributed according to this local equilibrium.  The existence and the uniqueness of the local equilibrium in $\Omega$ as well as the convergence toward this local equilibrium is made more precise in Section~\ref{sec.qsd} using the notion of quasi-stationary distribution.  To ensure the convergence towards the local equilibrium in~$\Omega$, a decorrelation step may be used before running the TAD algorithm, see step (M1)  in~\cite[Section 2.2]{aristoff-lelievre-14}.
 
 \noindent 
 As in the previous section, one denotes by~$(\Omega_i)_{i=1,\ldots,n}$ the  surrounding domains of~$\Omega$ (see Figure~\ref{fig:surrounding}), each of them corresponding to a macroscopic state of the system and, for $i\in \{1,\ldots,n\}$, $\{z_i\}={\rm argmin}_{\pa \Omega\cap \pa \Omega_i}f$. 
To sample the  next visited state among the~$\Omega_i$'s as well as the time~$T$ spent  in $\Omega$ for  the process~\eqref{eq.langevin}, the  TAD algorithm proceeds as follows. 
Let us introduce~$T_{sim}=0$ (which is the simulation  time) and $T_{stop}=+\infty$ (which is the stopping time), and iterate the following steps. 
\begin{sloppypar}
\begin{enumerate}
\item Let~$(Y_t)_{t\ge T_{sim}}$ be  the solution to the evolution equation~\eqref{eq.langevin} but for the temperature $h_{high}>h_{low}$,  starting from the local equilibrium in $\Omega$ at temperature $h_{high}$.  Let $(Y_t)_{t\ge T_{sim}}$ evolve until it leaves $\Omega$ and  denote by   $$T_{sim} +\tau$$ the first exit time  from~$\Omega$ for the process~$(Y_t)_{t\ge T_{sim}}$. Let $j\in\{1,\ldots,n\}$ be such that $Y_{T_{sim}+\tau}\in \pa \Omega_j\cap \pa \Omega$. Then, set  $T_{sim} =T_{sim} +\tau$.  If it is the first time  an exit   from $\Omega$ through $z_j$ for the process~$(Y_t)_{t\ge 0}$ is observed (else one goes directly to the next step), set $\tau_j(h_{high})=T_{sim}$ and extrapolate the time to~$\tau_j(h_{low})$ with the formula 
 \begin{equation}\label{tad-t}
\tau_j(h_{low})=\tau_j(h_{high}) \, e^{ 2\big(\frac{1}{h_{low}}-\frac{1}{h_{high}} \big ) (f(z_j)-f(x_0))  },
\end{equation}
 where we recall $x_0\in \Omega$ is the unique global minimum of $f$ in $\overline \Omega$. Then,  update the minimum exit time $\tau_{min}(h_{low})$ among the $\tau_j(h_{low})$'s which have been observed so far. Finally, compute a new time  $T_{stop}$  so that  there is a very small probability (say $\alpha \ll 1$) to observe an  exit event from~$\Omega$   at the temperature~$h_{high}$ which, using~\eqref{tad-t},  would change the value of $\tau_{min}(h_{low})$. We refer to~\cite{sorensen-voter-00} or~\cite{aristoff-lelievre-14}  for the computation of $T_{stop}$. 
\item  
 If  $T_{sim} \le T_{stop} $ then go back to the first step starting from the local  equilibrium in $\Omega$ at time $T_{sim}$, else go to the next step. 
 \item  Set $T=\tau_{min}(h_{low})$ and $Y= {\ell}$ where $\ell$ is such that $\tau_{\ell}(h_{low})=\tau_{min}(h_{low})$.  Finally, send~$X_{t_0+T}$ to $\Omega_\ell$ and evolve the process~\eqref{eq.langevin} with the new initial condition $X_{t_0+T}$.  
 \end{enumerate}
 \end{sloppypar}
 \begin{remark}\label{eq.reV}
 In~\cite{sorensen-voter-00},  when the process $(Y_t)_{t\ge T_{sim}}$ leaves $\Omega$, it is   reflected back in~$\Omega$ and it is then assumed that it reaches  instantaneously  the local equilibrium in~$\Omega$ at temperature $h_{high}$. 
  \end{remark}
  \begin{remark}
One can use a decorrelation step before running the TAD algorithm  and the sampling of    $Y_{T_{sim}}$ according to the local equilibrium  in $\Omega$ at the beginning of the step 1  to ensure that the  underlying Markov jump process is justified, see~\cite{aristoff-lelievre-14}. 
 \end{remark}
 \noindent
The extrapolation formula~\eqref{tad-t} which is at the heart of the  TAD  algorithm relies on the properties   of the underlying Markov jump process used to accelerate the exit event from a metastable state and  where  transition times are exponentially distributed with   parameters  computed with the Eyring-Kramers formula, see Remark~\ref{eq.equi-kmc} and Equation~\eqref{eq.ek}. 
In the algorithm TAD, it is indeed assumed that the exit event from  $\Omega$ can be modeled with a kinetic Monte Carlo method where the transition rates are computed with   the Eyring-Kramers law~\eqref{eq.ek}-\eqref{eq.prefactor}. Then, at high temperature, one checks that under this assumption, each     $\tau_ {i} (h_ {high})$ ($ i \in \{1, \ldots, n \}$) is an exponential law of parameter   $A_i\, e ^ {- \frac {2 } {h_ {high}}(f (z_ {i}) - f (x_0))}$ (see Remark~\ref{eq.equi-kmc}). The formula~\eqref{tad-t} allows to construct for all  $ i \in \{1, \ldots, n \}$, an exit time   $\tau_ {i} (h_ {low})$ which is  an exponential law of parameter $A_i\, e ^ {- \frac {2 } {h_ {low}}(f (z_ {i}) - f (x_0))}$. By considering the couple  $(\min_{ i\in \{1,\ldots,n\} }  \tau_{i}(h_ {low}), {\rm argmin}_{i\in \{1,\ldots,n\}}  \tau_{i}(h_ {low}))$, one has access to the exit event from~$\Omega$  (see Remark~\ref{eq.equi-kmc}).

\begin{remark}
There are other algorithms which use   the properties   of the underlying Markov jump process to accelerate the simulation of  the exit event from a metastable state, see for instance~\cite{voter-97} and~\cite{voter-98}. 
\end{remark}
Our objective is to  justify rigorously  that  a Markov jump process with transition rates   computed with the Eyring-Kramers formula~\eqref{eq.ek}  can be used to model   the exit event from a metastable domain~$\Omega$ for the overdamped Langevin process~\eqref{eq.langevin}.   Before, let us  recall mathematical contributions  on  the  exit event from a  domain and   on the Eyring-Kramers formula~\eqref{eq.ek}.

 \subsection{Mathematical   literature on the  exit event from a  domain and  on the Eyring-Kramers formulas} 
In the mathematical  literature, there are  mainly two  approaches to the study of the asymptotic behaviour of the exit event from a domain when $h\to 0$: the global approaches and the local approaches.   

\subsubsection{Global approaches.}
 
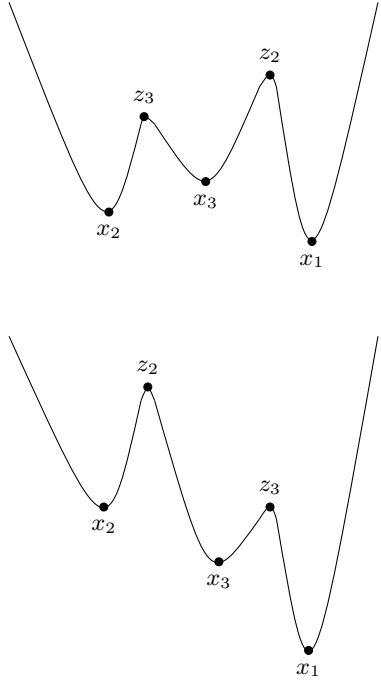
\begin{figure}[h!]
\begin{center}
\begin{tikzpicture}[scale=0.45]
\tikzstyle{vertex}=[draw,circle,fill=black,minimum size=3pt,inner sep=0pt]
\draw (-5.9,5) ..controls  (-3,-2.5).. (-2,1.4)  ;
\draw (-1.6,1.4) ..controls  (0,-1).. (1.5,2.5)  ;
\draw (5,5) ..controls  (3,-3.9).. (2,2.5)  ;
\draw (-2,1.4) ..controls  (-1.91,1.7).. (-1.6,1.4)   ;
\draw (2,2.5) ..controls  (1.83,3).. (1.5,2.5)   ;
\draw (3.05,-2.07) node[vertex,label=south: {$x_1$}](v){};
\draw (-2.95,-1.2) node[vertex,label=south: {$x_2$}](v){};
\draw (-0.09,-0.3)  node[vertex,label=south: {$x_3$}](v){};
\draw (1.81,2.85) node[vertex,label=north: {$z_2$}](v){};
\draw (-1.91,1.62) node[vertex,label=north: {$z_3$}](v){};
\end{tikzpicture}

\hspace{1cm}

\begin{tikzpicture}[scale=0.45]
\tikzstyle{vertex}=[draw,circle,fill=black,minimum size=3pt,inner sep=0pt]
\draw (-5.9,5.9) ..controls  (-3,-0.5).. (-2,4)  ;
\draw (-1.6,4) ..controls  (0,-1.6).. (1.5,0.5)  ;
\draw (5,5.9) ..controls  (3,-5.4).. (2,0.5)  ;
\draw (-2,4) ..controls  (-1.8,4.5).. (-1.6,4)   ;
\draw (2,0.5) ..controls  (1.83,1).. (1.5,0.5)   ;
\draw (2.95,-3.39) node[vertex,label=south: {$x_1$}](v){};
\draw (-3.1,0.85) node[vertex,label=south: {$x_2$}](v){};
\draw (0.3,-0.77)  node[vertex,label=south: {$x_3$}](v){};
\draw (1.81,0.85) node[vertex,label=north: {$z_3$}](v){};
\draw (-1.8,4.4) node[vertex,label=north: {$z_2$}](v){};
\end{tikzpicture}

\caption{Examples of two labelings  of the local minima $\{x_1,x_2,x_3\}$ of $f$ in dimension one.   }
 \label{fig:meta-order}
 \end{center}
\end{figure}

The global approaches   study the asymptotic behaviours in the limit $h\to 0$ of the eigenvalues of the infinitesimal generator  
\begin{equation}\label{eq.ig}
  L^{ (0)}_{f,h} =-\frac h2 \Delta +\nabla f\cdot \nabla 
\end{equation}
 of the diffusion~\eqref{eq.langevin} on $\mathbb R^d$.  Let us give for example a result obtained in~\cite{bovier-eckhoff-gayrard-klein-04,bovier-gayrard-klein-05}. 
To this end, let us assume that the  potential $f:\mathbb R^d\to \mathbb R$ is a Morse function,   has $m$ local minima  $\{x_1,\ldots,x_m\}$   and that for $h$ small enough $\int_{\mathbb R^d}e^{-\frac 2h f}<+\infty$. 
Let us recall that $\phi:\mathbb R^d\to \mathbb R$ is a Morse function if all its critical points are non degenerate. For a Morse function  $\phi:\mathbb R^d\to \mathbb R$, we say that $x$ is a saddle point of $\phi$ if $x$ is a critical point of~$\phi$ such that the Hessian matrix of $\phi$  at $x$ has exactly one negative eigenvalue (i.e. $x$ is a critical point of $\phi$ of index $1$).
Then, from~\cite{HeSj4}, the operator~$L^{ (0)}_{f,h} $ has exactly  $m$ exponentially small eigenvalues  $\{\lambda_1,\lambda_2,\ldots,\lambda_m\}$ when $h\to 0$ with $\lambda_1=0<\lambda_2\le \ldots\le \lambda_m$ (i.e., when $h\to 0$,  for all $i\in\{1,\ldots,m\}$, $\lambda_i=O( e^{-\frac ch})$ for some $c>0$ independent of $h$).   
Moreover, sharp asymptotic estimates can be derived for the eigenvalues~$\{\lambda_2,\ldots,\lambda_m\}$.  In~\cite{bovier-eckhoff-gayrard-klein-04,bovier-gayrard-klein-05}, the following  results  are obtained.  Let us assume that  $\{x_1\}={\rm argmin }_{\mathbb R^d} f$.  For $k\in \{2,\ldots,m\}$ and $B_k=\{x\in \{x_1,\ldots,x_m\}\setminus\{x_k\}, f(x)\le f(x_k)\} $ (i.e. $B_k$ is the set of local minima of $f$ which are lower in energy than $x_k$), one denotes by~$\mathcal P(x_k,B_k)$  the set of curves    $\gamma\in C^0([0,1], \mathbb R^d)$ such that  $\gamma(0)=x_k$ and $\gamma(1)\in B_k$. Let us finally assume that:
\begin{enumerate}
\item For all $k\in \{2,\ldots,m\}$, there exists a unique saddle point~$z_k$  (i.e. a critical point of $f$ of index $1$) such that  $f(z_k)=\inf_{\gamma\in \mathcal P(x_k,B_k)}\sup_{t\in [0,1]} f(\gamma(t))$.   
\item \begin{sloppypar}   The values  $\big (f(z_k)-f(x_k)\big)_{k\in\{2,\ldots,m\}}$ are all distinct. \end{sloppypar}
\end{enumerate}
\begin{sloppypar}
\noindent
These assumptions imply that the map $x_k\in \{x_2,\ldots,x_m\}  \mapsto z_k$ is injective. 
The set $\{x_2,\ldots,x_m\}$ is then labeled such that the sequence  $\big (f(z_k)-f(x_k)\big)_{k\in\{2,\ldots,m\}}$ is strictly decreasing. 
The previous assumptions  also imply   the existence of a  cascade of events, which occur with  different timescales,   to go   from one local minimum $x_k$ of $f$ to the global minimum~$x_1$ of~$f$ in $\mathbb R^d$, see for instance Figure~\ref{fig:meta-order}. 
Then,  one has for $k\in \{2,\ldots,m\}$,  in the limit $h\to 0$:
\end{sloppypar}
\begin{equation}\label{eq.lamdakk}
\lambda_k= \frac{\vert \lambda(z_k)\vert}{ 2\pi}\frac{    \sqrt{ {\rm det\,  Hess}f (x_k) }    }{  \vert  \sqrt{ {\rm det \,Hess}f(z_k)  \vert } } e^{-\frac 2h (f(z_k)-f(x_k))} (1+o(1)),
\end{equation}
where   $\lambda(z_k)$ is the negative eigenvalue of the Hessian matrix of~$f$ at~$z_k$.  In the articles~\cite{bovier-eckhoff-gayrard-klein-04,bovier-gayrard-klein-05}, using a potential-theoretic approach, the sharp equivalent~\eqref{eq.lamdakk} is obtained and 
each of the eigenvalues~$\lambda_k$ (for $k\in \{2,\ldots,m\}$) is shown to be the inverse of the average time it takes for the process~\eqref{eq.langevin} to go from~$x_k$ to  $B_k$.  We also refer to~\cite{eckhoff-05} for similar results. 
In~\cite{helffer-klein-nier-04},  another proof of~\eqref{eq.lamdakk} is given using tools from semi-classical analysis. Let us also mention~\cite{michel2017small} for a generalization of the results obtained in~\cite{helffer-klein-nier-04}. Notice that the  results presented above  do not   provide any information  concerning the average  time it takes for the process~\eqref{eq.langevin} to go from the global minimum of~$f$ to   a  local minimum of $f$ when $h\to 0$. One also refers to~\cite{landim2017dirichlet} for generalization of~\cite{bovier-eckhoff-gayrard-klein-04,bovier-gayrard-klein-05} for a class of non reversible processes when $f$ has two local minima, and to~\cite{miclo-95,holley-kusuoka-stroock-89,davies-82a,davies-82b,davies-82c} for  related results.   
 
\begin{remark}The global approaches have been used in~\cite{schuette-sarich-13,schuette-98} to construct a Markovian dynamics by projecting  the infinitesimal generator $L_{f,h}^{(0)}$ of the diffusion~\eqref{eq.langevin} with a Galerkin method onto the vector space associated with the  $m$ small eigenvalues~$\{\lambda_1,\ldots,\lambda_m\}$. This projection leads to  a very good approximation of~$L_{f,h}^{(0)}$  in the limit~$h\to 0$.  The question is then how to relate the transition
events (or the trajectories) of the obtained Markov
process to the exit events (or the trajectories) of
the original one.
\end{remark}

\subsubsection{Local approaches.} The local approaches consist in studying the asymptotic behaviour when $h\to 0$ of the exit event $(\tau_\Omega, X_{\tau_\Omega})$ from a domain~$\Omega\subset \mathbb R^d$, where $\tau_\Omega:=\inf\{t\ge 0, X_t\notin \Omega \}$.  \medskip

\noindent
One of the most well-known approaches  is the large deviation theory developed by  Freidlin and
Wentzell in the 1970s. We refer to the book~\cite{freidlin-wentzell-84} which  summarizes their main contributions. This theory is based on the study of small  pieces of the trajectories of the process defined with a suitable increasing sequence of stopping times. The rate function is fundamental in this theory: it quantifies  the cost  of deviating from a  deterministic trajectory when $h\to 0$. The  rate functional was first introduced by~Schilder~\cite{schilder1966some} for a Brownian motion.  Some typical results from~\cite{freidlin-wentzell-84} (see Theorem~2.1,  Theorem~4.1, and Theorem~5.1 there) are the following. Let~$\Omega$  be a $C^{\infty}$ open and  connected bounded subset of $\mathbb R^d$. Let us assume that   $\pa_nf>0$ on $\pa\Omega$ (where $\partial_n$ is the outward normal derivative to $\Omega$) and  that $f$ has a unique non degenerate critical point  $x_0$ in $\Omega$   such that   $f(x_0)=\min_{\overline \Omega} f$. Then, for all~$x\in \Omega$:
$$\lim_{h\to 0} h\ln \mathbb E_x\big [\tau_\Omega\big ]=2\big (\inf_{\pa\Omega}f -f(x_0)\big ).$$
 The notation  $\mathbb E_x$ stands for the expectation given the fact that $X_0=x$. 
Moreover, let~$x\in \Omega$ such that~$f(x)< \inf_{\pa\Omega}f$. Then, for any $\gamma>0$ and  $\delta_0>0$, there exist~$\delta\in (0,\delta_0]$  and $h_0>0$ such that  for all $h\in (0,h_0)$ and for all~$y\in \pa \Omega$:
$$   e^{-\frac 2h (f(y)-\inf_{\pa\Omega}f )}e^{-\frac \gamma h}\le \mathbb P_x\big [\vert X_{\tau_\Omega}-y\vert <\delta \big ]\le e^{-\frac 2h (f(y)-\inf_{\pa\Omega}f )}e^{\frac \gamma h}.$$
The notation  $\mathbb P_x$ stands for the probability given the fact that $X_0=x$. 
Lastly, if the infimum of  $f$ on $\pa \Omega$ is attained at one single point $y_0\in \pa \Omega$, then for all $\delta>0$:
$$\lim_{h\to 0} \mathbb P_x\big [\vert X_{\tau_\Omega}-y_0\vert <\delta \big ]=1.$$
  A  result due to Day~\cite{day-83}  (see also \cite{mathieu-95,mathieu-94}) concerning the law of $\tau_\Omega$ is the following. When~$h\to 0$, the exit time~$\tau_\Omega$ converges in law to an exponentially distributed random variable  and for all~$x\in \Omega$
$$ \lim_{h\to 0} \lambda_h \mathbb E_x\big [\tau_\Omega\big ] =1,$$
where $\lambda_h$ is the principal eigenvalue of the infinitesimal generator of the diffusion~\eqref{eq.langevin} associated with Dirichlet boundary conditions on $\pa \Omega$ (see Proposition~\ref{fried-bis} below). 
The interest of this approach is that it can be applied to very general dynamics.
However, when it is used to prove that the Eyring-Kramers formulas~\eqref{eq.ek} can be used to
study the exit distribution from $\Omega$, it only provides the exponential rates (not
the prefactor $A_i$ in~\eqref{eq.ek}) and does not give error bounds when $h \to 0$.
\medskip

\noindent
There are also approaches which are based on techniques developed for partial differential equations.   In~\cite{MS77,MaSc}, using formal computations, when  $\pa_nf>0$ on $\pa\Omega$ and   $f$ has a unique non degenerate critical point  $x_0$ in $\Omega$   such that   $f(x_0)=\min_{\overline \Omega} f$,   the following formula is derived: for any~$F\in C^{\infty}(\pa \Omega,\mathbb R)$ and $x\in \Omega$, one has when $h\to 0$:
\begin{equation}\label{eq.MS}
 \mathbb E_x\big [F\big (X_{\tau_\Omega}\big) \big ]=\frac{\int_{\pa \Omega} F(z) \partial_nf(z)\, e^{-\frac{2}{h} f(z)}\, dz  }{\int_{\partial \Omega} \partial_nf \, e^{-\frac{2}{h} f } d\sigma }  +o(1).
 \end{equation}
The formal asymptotic estimate~\eqref{eq.MS} implies that the law of  $ X_{\tau_\Omega}$ concentrates on points where $f$ attains its minimum on~$\pa \Omega$. Moreover,  an asymptotic equivalent of $\mathbb E_x\big [\tau_\Omega\big ] $ when $h\to 0$ is also formulated in~\cite{schuss90} through formal computations.   These results are obtained injecting formal asymptotic expansions in powers of $h$  in  the  partial differential equations satisfied by $x\in \Omega\mapsto  \mathbb E_x\big [F\big (X_{\tau_\Omega}\big) \big ]$ and $x\in \Omega\mapsto  \mathbb E_x\big [ \tau_\Omega \big ]$. 
We also refer to~\cite{MaSc}, where using formal computations,   asymptotic formulas  are obtained  concerning  both  the concentration of the law of $X_{\tau_\Omega}$ on ${\rm argmin}_{\pa \Omega} f$ and   $ \mathbb E_x\big [ \tau_\Omega \big ]$  when $\Omega$ is the union of basins of attraction of the dynamics $\frac{d}{dt}\gamma(t)=-\nabla f(\gamma(t))$. 
When  $\pa_nf>0$ on $\pa\Omega$ and   $f$ has a unique non degenerate critical point  $x_0$ in $\Omega$   such that   $f(x_0)=\min_{\overline \Omega} f$, the formula~\eqref{eq.MS}  is proved rigorously by Kamin in~\cite{kamin1979elliptic}, and is extended to a non reversible  diffusion  process $(Y_t)_{t\ge 0}$ solution to $dY_t=b(Y_t)\, dt +\sqrt{h}\, dB_t$ in~\cite{kamin-78,perthame-90,day1984a,day1987r} when~$\Omega$ contains  one attractor of the dynamics $\frac{d}{dt}\gamma(t)=b(\gamma(t))$ and $b(x)\cdot n<0$ for all~$x\in \pa \Omega$. However, the results~\cite{kamin1979elliptic,kamin-78,perthame-90,day1984a,day1987r}  do not provide any information on  the probability to leave~$\Omega$ through a point which is not a global minimum of~$f$ on~$\pa \Omega$. \medskip

\noindent
 Finally, let us mention~\cite{devinatz-friedman-78a,devinatz-friedman-78b,holley-kusuoka-stroock-89,mathieu-95,mathieu-94, helffer-nier-06, le-peutrec-10} for a study of the asymptotic behaviour in the limit $h\to 0$ of  $\lambda_h$ and $u_h$ (see Proposition~\ref{fried-bis} below).  The reader can also refer to~\cite{day-99} for a review of the different techniques used to  study the asymptotic behaviour of  $ X_{\tau_\Omega}$ when $h\to 0$ and to~\cite{berglund-13} for a review of the different techniques used to  study the asymptotic behaviour of  $\tau_\Omega$ when $h\to 0$.  
 \begin{remark}
   \begin{sloppypar}
 Some authors proved the convergence to a Markov jump process in some specific geometric settings and after a rescaling in
time. We refer to~\cite{kipnis-newman-85}  for a one-dimensional diffusion in a double well and~\cite{galves-olivieri-vares-87,mathieu-95} for a study in higher dimension. In~\cite{sugiura-95}, assuming that all the saddle points of $f$ are at the same height, it is proved that a   suitable  rescaling of the time leads to a convergence of the diffusion process to   a Markov jump process  between the global minima of~$f$.
  \end{sloppypar}
\end{remark}
The results presented in this work (see~\cite{di-gesu-lelievre-le-peutrec-nectoux-16,DLLN}) follow   a local approach. The quasi-stationary distribution of the process~\eqref{eq.langevin} on  $\Omega$ is the cornerstone of the analysis. They state that, under some geometric assumptions, the Eyring-Kramers
formulas (with prefactors) can be used to model the exit event from a
metastable state, and provide explicit error bounds.

  \subsection{Quasi-stationary distribution and  transition rates} 
  \label{sec.qsd}
  \subsubsection{Local equilibrium.}
  Let~$\Omega$ be a $C^\infty$  open bounded connected subset of $\mathbb R^d$ and $f\in C^\infty(\overline \Omega, \mathbb R)$. Let us recall that $
\tau_\Omega:=\inf\{t\ge 0, X_t\notin \Omega \}$ denotes  
  the first exit time from~$\Omega$.  The quasi-stationary distribution of the process~\eqref{eq.langevin} on  $\Omega$ is defined as follows. 
  \begin{definition}\label{de.qsd}
A probability measure  $\nu_h$ on $\Omega$ is a quasi-stationary distribution of the process~\eqref{eq.langevin} on  $\Omega$ if for all $t> 0$ and any measurable set $A\subset\Omega$, 
 $$\mathbb P_{\nu_h}\big[X_t\in A\big | t<\tau_\Omega   \big]=\nu_h(A).$$
\end{definition}  
  The notation  $\mathbb P_\mu$ stands for the probability given the fact that the process~\eqref{eq.langevin} is initially distributed according to~$\mu$ i.e. $X_0\sim \mu$. The next   proposition~\cite{le-bris-lelievre-luskin-perez-12,champagnat2015quasi} shows that the law of the process~\eqref{eq.langevin} at time $t$ conditioned not to leave $\Omega$ on the interval $(0,t)$  converges to the quasi-stationary distribution. 
\begin{proposition}\label{pr.conv}
Let~$\Omega$ be a $C^\infty$ open connected and bounded subset of $\mathbb R^d$ and $f\in C^\infty(\overline \Omega, \mathbb R)$.  Then, there exist a unique probability measure  $\nu_h$ on $\Omega$  and  $c>0$ such that for any probability measure $\mu$ on $\Omega$, there exist $C(\mu)>0$  and $t(\mu)>0$  such that for all $t\ge t(\mu)$  and all measurable set $A\subset\Omega$:
\begin{equation}\label{eq.con-c}
\big \vert \mathbb P_\mu\big[X_t\in A\big | t<\tau_\Omega   \big] -\nu_h(A)      \big \vert \le  C(\mu)e^{-ct}.
\end{equation}
Moreover, $\nu_h$ is the unique quasi-stationary distribution of the process~\eqref{eq.langevin} on  $\Omega$. 
\end{proposition}
Proposition~\ref{pr.conv} indicates that the quasi-stationary distribution~$\nu_h$ can be seen as a local equilibrium of the process~\eqref{eq.langevin} in~$\Omega$. 
\medskip

\noindent
The quasi-stationary distribution~$\nu_h$ can be expressed with the principal eigenfunction of the infinitesimal generator $L^{ (0)}_{f,h}$ (see~\eqref{eq.ig}) of the diffusion~\eqref{eq.langevin} associated with Dirichlet boundary conditions on~$\pa \Omega$.  To this end, let us introduce the following Hilbert spaces $L^2_w(\Omega)=\big \{u:\Omega \to \mathbb R,  \, \int_\Omega u^2
 e^{-\frac{2}{h} f }   < \infty\big\}$  
and for $q\in \{1,2\}$,
\begin{equation}\label{eq.H1w}
H^{q}_w(\Omega)= \big\{u \in L^2_w(\Omega), \,  \forall \alpha\in \mathbb N^d,\,  \vert \alpha\vert \leq q,  \, \pa_\alpha u \in L^2_w(\Omega)\big  \}.
\end{equation}
The subscript $w$ in the notation $L^2_w(\Omega)$ and $H^{q}_w(\Omega)$ refers to the
fact that the weight  function~$x\in \Omega \mapsto e^{-\frac{2}{h} f(x)} $ appears in the inner product. 
Moreover, let us denote by $H^{1}_{0,w}(\Omega)=\{u\in H^{1}_{w}(\Omega), \, u=0 \text{ on } \pa \Omega\}$. 
Let us recall the following result~\cite{le-bris-lelievre-luskin-perez-12}. 
\begin{proposition}\label{fried-bis}
 Let~$\Omega$ be a $C^\infty$ open connected and bounded subset of~$\mathbb R^d$ and~$f\in C^\infty(\overline \Omega, \mathbb R)$. Then, the operator $L^{ (0)}_{f,h}$ with domain $H^{1}_{0,w}(\Omega)\cap H^{2}_w(\Omega)$ on~$L^2_w(\Omega)$, which is denoted by $L^{D, (0)}_{f,h}$, is self-adjoint, positive and has compact resolvent. Furthermore, the smallest eigenvalue $\lambda_h$ of $L^{D, (0)}_{f,h}$ is non degenerate and any eigenfunction associated with $\lambda_h$ has a sign on $\Omega$. 
  \end{proposition}
  In the following, one denotes by~$u_h$ an  eigenfunction associated with $\lambda_h$. The smallest eigenvalue $\lambda_h$ of $L^{D, (0)}_{f,h}$ is called the principal eigenvalue   of~$L^{D, (0)}_{f,h}$ and~$u_h$ a principal eigenfunction of~$L^{D, (0)}_{f,h}$. 
   Without loss of generality, one assumes that
  \begin{equation}\label{norma}
u_h > 0 \text{ on } \Omega \ \text{ and }  \int_{\Omega} u_h^2\, e^{-\frac{2}{h}f}=1.
\end{equation}
Then, the quasi-stationary distribution~$\nu_h$ of the process~\eqref{eq.langevin} in~$\Omega$ is given by (see~\cite{le-bris-lelievre-luskin-perez-12}): 
\begin{equation}\label{exp.qsd}
\nu_h(dx)=\frac{u_h(x)\ e^{-\frac{2}{h}f(x)} }{\int_{\Omega} u_h\, e^{-\frac{2}{h}f}}\, dx. 
\end{equation}
Moreover,  the following result shows that when $X_0\sim \nu_h$, the law of the exit event $(\tau_\Omega, X_{\tau_\Omega})$  is explicitly known  in terms of $\lambda_h$ and $u_h$ (see~\cite{le-bris-lelievre-luskin-perez-12}). 

  \begin{proposition}\label{indep1-1}
Let us assume that $X_0\sim \nu_h$, where $\nu_h$ is the quasi-stationary distribution of the process~\eqref{eq.langevin} in~$\Omega$. Then,
$\tau_{\Omega}$ and $X_{\tau_{\Omega}}$ are independent. Moreover, $\tau_{\Omega}$ is exponentially distributed  with parameter  $\lambda_h$ and for any open set~$\Sigma\subset \pa \Omega$, one has:
\begin{equation}\label{expression-P}
\mathbb P_{\nu_{h}}\big [ X_{\tau_{\Omega}}\in \Sigma \big]=-\frac {h}{2\lambda_h} \frac{\displaystyle \int_{\Sigma}   \partial_n
  u_h(z) e^{-\frac{2}{h} f(z)} \sigma(dz) }{\displaystyle \int_\Omega u_h e^{-\frac{2}{h} f}},
    \end{equation}
    where $\sigma(dz)$ is the Lebesgue measure on $\pa \Omega$. 
\end{proposition}
  
  \subsubsection{Approximation of the exit event with a Markov jump process.} Let us now provide justifications to the use of  a Markov jump process with transition rates   computed with the Eyring-Kramers formula~\eqref{eq.ek} to model   the exit event from a metastable domain~$\Omega$. 
 In view of~\eqref{eq.con-c}, one can be more precise on the definition of  the metastability of a domain~$\Omega$ given in Section~\ref{sec11}.  For a probability measure $\mu$ supported in~$\Omega$, the domain~$\Omega$ is said to be metastable  if, when $X_0\sim \mu$, the convergence to the quasi-stationary distribution~$\nu_h$ in~\eqref{eq.langevin} is much quicker than the  exit  from~$\Omega$.
Since the process~\eqref{eq.langevin} is a Markov process, it is then  relevant to study the exit event from~$\Omega$ starting from the quasi-stationary distribution~$\nu_h$, i.e. $X_0\sim\nu_h$. As a consequence of Proposition~\ref{indep1-1}, the exit time is exponentially distributed and is independent of the next visited state. These two properties are  the fundamental  features of   kinetic Monte Carlo methods, see indeed~\eqref{eq.Te} and~\eqref{eq.Ye}.  It thus remains to  prove that the transition rates can be computed with the Eyring-Kramers  formula~\eqref{eq.ek}. For that purpose, let us first give an expression of the transition rates. 
Recall that~$(\Omega_i)_{i=1,\ldots,n}$ denotes the  surrounding domains of~$\Omega$ (see Figure~\ref{fig:surrounding}). For $i\in\{1,\ldots,n\}$, we define  the transition rate to go from $\Omega$ to $\Omega_i$ as follows:
 \begin{equation}\label{kijl}
k_{i}^L:=\frac{1}{\mathbb E_{\nu_{h}}\big [ \tau_{\Omega}\big]} \mathbb P_{\nu_{h}}\big [ X_{\tau_{\Omega}}\in \pa \Omega\cap \pa \Omega_i \big],
\end{equation}
where we recall, $\nu_h$ is the  quasi-stationary distribution of the process~\eqref{eq.langevin} in~$\Omega$. The superscript $L$ in~\eqref{kijl} indicates that the microscopic evolution of the system is governed by  the overdamped Langevin process~\eqref{eq.langevin}. Notice that, using  Proposition~\ref{indep1-1}, it holds for all $i\in \{1,\ldots,n\}$: 
$$
\mathbb P_{\nu_{h}}\big [ X_{\tau_{\Omega}}\in \pa \Omega\cap \pa \Omega_i \big] \ =\ \frac{ k_{i}^L}{\sum_{\ell=1}^n k_{\ell}^L}.
$$
Thus, the expressions~\eqref{kijl} are compatible with the use of a kinetic Monte Carlo algorithm, see~\eqref{eq.Te} and~ \eqref{eq.Ylaw}. Indeed, starting from the quasi-stationary distribution~$\nu_h$, the
exit event from $\Omega$ can be exactly modeled using
the rates~\eqref{kijl}: the exit time is   exponentially
distributed with parameter $\sum_{\ell=1}^n k_{\ell}^L$, independent of the
exit point, and the exit point is in $\partial \Omega_i \cap
\partial \Omega$ with probability $k^L_i / \sum_{\ell=1}^n k_{\ell}^L$. The remaining question is thus following: does the transition rate~\eqref{kijl} satisfy the Eyring-Kramers law~\eqref{eq.ek} in the limit~$h\to 0$?  

Notice that, using Proposition~\ref{indep1-1}, for $i\in \{1,\ldots,n\}$, the transition rate defined by~\eqref{kijl} writes:
\begin{equation}\label{expression-kijl}
k_i^L=-\frac h2 \frac{\displaystyle \int_{\pa \Omega\cap \pa \Omega_i}   \partial_n
  u_h(z) \, e^{-\frac{2}{h} f(z)} \sigma(dz) }{\displaystyle \int_\Omega u_h \, e^{-\frac{2}{h} f}},
  \end{equation}
where we recall, $u_h$ is the eigenfunction associated with the principal eigenvalue~$\lambda_h$  of~$L^{D, (0)}_{f,h}$.  
\medskip

\begin{sloppypar}
\noindent 
The remainder of this work is dedicated to the presentation of recent results in~\cite{di-gesu-lelievre-le-peutrec-nectoux-16},~\cite{DLLN} and~\cite{BN2017}  which aim at studying  the  asymptotic behaviour of   the exit event~$(\tau_\Omega, X_{\tau_\Omega})$ from a metastable domain~$\Omega$ in the limit $h\to 0$. 
In particular,  the results give a sharp  asymptotic formula of  the transition rates~\eqref{expression-kijl} when~$h~\to~0$. 
\end{sloppypar}

\begin{remark}\label{re.facteur2}
If one wants to recover the expression of the prefactor~\eqref{eq.prefactor}, one has to multiply by $\frac 12$ the expression~\eqref{kijl}. This can be explained as follows.  Once the process~\eqref{eq.langevin} reaches $ \pa \Omega\cap \pa \Omega_i $, it has, in the limit $h\to 0$,  a one-half  probability to come back in  $\Omega$  and  a one-half  probability   to go in $\Omega_i$. If $z_i$ is a non degenerate saddle point of $f$, this result is not difficult to prove in dimension $1$. Indeed, it is proved in~\cite[Section A.1.2.2]{BN2017},  that when reaching $ \pa \Omega\cap \pa \Omega_i $, the probability that the process~\eqref{eq.langevin}    goes in  $\Omega_i$ is   $\frac 12+O(h)$  in the limit $h\to 0$. To extend this result to  higher dimensions,  one can use a suitable set of coordinates around~$z_i$.  
\end{remark}

\section{Main results on the exit event} 
\label{sec.mainresults}

In all this section,~$\Omega\subset \mathbb R^d$ is~$C^\infty$ open, bounded and connected, and~$f\in C^\infty(\overline \Omega,\mathbb R)$
\footnote{Actually, all the results presented in this section are proved in~\cite{DLLN} and~\cite{di-gesu-lelievre-le-peutrec-nectoux-16} in the   more general setting: $\overline\Omega=\Omega\cup \pa \Omega$ is a $C^{\infty}$ oriented compact and  connected Riemannian manifold of dimension $d$ with boundary~$\partial \Omega$.}. 
The purpose of this section is to present recent results obtained in~\cite{DLLN} and~\cite{di-gesu-lelievre-le-peutrec-nectoux-16}. Both~\cite{DLLN} and~\cite{di-gesu-lelievre-le-peutrec-nectoux-16} are mainly concerned with studying  the asymptotic behaviour when $h\to 0$ of  the exit law of a domain~$\Omega$ of the process~\eqref{eq.langevin}.  In~\cite{DLLN}, when~$\Omega$ only contains  one local minimum of $f$ and $\pa_n f> 0$ on~$\pa \Omega$, we obtain sharp asymptotic equivalents when $h\to 0$ of the probability that the process~\eqref{eq.langevin} leaves~$\Omega$  through a subset $\Sigma$ of  $\pa \Omega$ starting from the quasi-stationary distribution or from a deterministic initial condition in~$\Omega$. Then, these asymptotic equivalents are used to compute the asymptotic behaviour of the transition rates~\eqref{kijl}. In~\cite{di-gesu-lelievre-le-peutrec-nectoux-16}, we explicit a more general setting than the one considered in~\cite{DLLN} where we identify  the most probable places of exit of $\Omega$ as well as their relative probabilities starting from the quasi-stationary distribution  or    deterministic initial conditions in~$\Omega$. More precisely, we consider in~\cite{di-gesu-lelievre-le-peutrec-nectoux-16} the case when~$\Omega$ contains several local minima of~$f$  and~$\vert \nabla f\vert\neq 0$ on $\pa \Omega$.  

\subsection{Sharp asymptotic estimates on the exit event from a domain} 
\label{sec.m1}
In this section, we present the results of~\cite{DLLN} which give sharp asymptotic estimates on the law of $X_{\tau_\Omega}$ and on the expectation of $\tau_\Omega$ when $h\to 0$. 
These results   give in particular  the asymptotic estimates of the transition rates $(k_j^L)_{j=1,\ldots,n}$ defined in~\eqref{kijl}. 

\subsubsection{Geometric setting.}

 Let us   give the  geometric setting which is considered in this section:
\begin{itemize}
\item \textbf{[H1]} The function $f:\overline \Omega\to \mathbb R$ and the restriction of $f$ to $\Omega$, denoted by~$f|_{\partial \Omega}$, are Morse functions. Moreover, $\vert \nabla f\vert (x)\neq 0$ for all $x\in \partial \Omega$.
\item \textbf{[H2]} The function $f$ has a unique global minimum $x_0$ in $\overline \Omega$ and 
$$\min_{\partial \Omega }f>\min_{\overline \Omega }f=\min_{ \Omega }f=f(x_0).$$
The point $x_0$ is the unique critical point of $f$ in $\overline \Omega$. The  function $f|_{\partial \Omega}$ has exactly $n\ge 1$ local minima which are denoted by $(z_i)_{i=1,\ldots,n}$. They are ordered such that $$f(z_1)\le \ldots\le f(z_n).$$
\item \textbf{[H3]} $\partial_nf(x)>0$ for all $x\in \partial \Omega$.
\end{itemize}

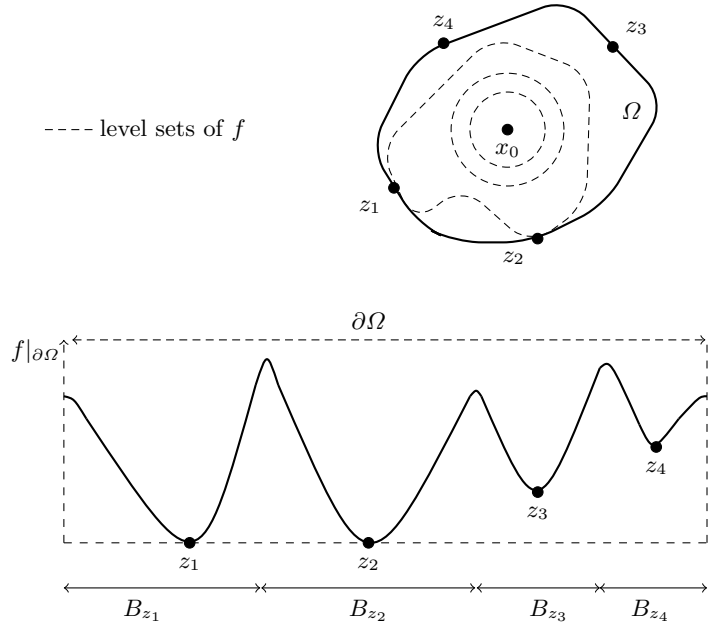
\begin{figure}[h!]
\begin{center}
\begin{tikzpicture}[scale=0.5]
\tikzstyle{vertex}=[draw,circle,fill=black,minimum size=4pt,inner sep=0pt]
\tikzstyle{ball}=[circle, dashed, minimum size=1cm, draw]
\tikzstyle{point}=[circle, fill, minimum size=.01cm, draw]
\draw [thick, rounded corners=10pt] (1,0.5) -- (-0.25,2.5) -- (1,5) -- (5,6.5) -- (7.6,3.75) -- (6,1) -- (4,0) -- (2,0) --cycle;
 
\draw [ densely dashed, rounded corners=10pt] (0.9,0.5) -- (-0,2.5) -- (3,5.5) -- (5.6,4.5) --(5.5,1) -- (4,-0.1) -- (2.2,1.5) --cycle;

\draw [densely dashed] (3.4,3) circle (1);
\draw [densely dashed] (3.4,3) circle (1.5);

\draw [densely dashed] (-8.9,3) -- (-7.6,3);
     \draw  (-5.5,3) node[]{$\text{level sets of $f$}$};

     \draw  (6.7,3.5) node[]{$\Omega$};
   
\draw (1.7,5.3) node[vertex,label=north: {$z_4$}](v){};
\draw (3.4,3) node[vertex,label=south: {$x_0$}](v){};
\draw (4.2,0.1) node[vertex,label=south west: {$z_2$}](v){};
\draw (0.38,1.45) node[vertex,label=south west: {$z_1$}](v){};
\draw (6.2,5.2) node[vertex,label=north east: {$z_3$}](v){};
\end{tikzpicture}

\hspace{0.7cm}

\begin{tikzpicture}[scale=0.75]
\tikzstyle{vertex}=[draw,circle,fill=black,minimum size=4pt,inner sep=0pt]
\tikzstyle{ball}=[circle, dashed, minimum size=1cm, draw]
\tikzstyle{point}=[circle, fill, minimum size=.01cm, draw]

\draw [dashed] (-5.4,-0.6)--(6,-0.6);
\draw [dashed,->] (-5.4,-0.6)--(-5.4,3);
\draw [dashed] (6,-0.6)--(6,3);
\draw [densely dashed,<->] (-5.25,3)--(5.97,3);
 \draw (0,3.34) node[]{$\pa \Omega$};
 \draw (-5.9,2.8) node[]{$f|_{\pa \Omega}$};
 \draw[thick] (-5.4,2) ..controls  (-5.22,1.96).. (-5,1.6)  ;
\draw[thick] (-5,1.6) ..controls  (-3,-1.4).. (-2,2.2)  ;
\draw[thick] (-1.6,2.2) ..controls  (0,-1.5).. (1.8,2)  ;
\draw[thick] (-2,2.2) ..controls  (-1.8,2.81).. (-1.6,2.2)   ;
\draw[thick] (2,2) ..controls  (1.91,2.13).. (1.8,2)   ;
\draw[thick] (4.1,2.5) ..controls  (3,-0.3).. (2,2)  ;
\draw[thick]  (4.1,2.5) ..controls  (4.28,2.66).. (4.8,1.5)  ;
\draw[thick]  (4.8,1.5)  ..controls  (5.04,1).. (5.5,1.6)  ;
\draw[thick]  (5.5,1.6)   ..controls  (5.87,2).. (6,2)  ;

\draw (5.1,1.1) node[vertex,label=south: {$z_4$}](v){};
\draw (3,0.3) node[vertex,label=south: {$z_3$}](v){};
\draw (-3.17,-0.6) node[vertex,label=south: {$z_1$}](v){};
\draw (0,-0.6)  node[vertex,label=south: {$z_2$}](v){};

\draw[<->]  (-5.4,-1.4) -- (-1.92,-1.4) ;
\draw[<->]  (-1.9,-1.4) -- (1.89,-1.4) ;
\draw[<->]   (1.92,-1.4) -- (4.1,-1.4) ;
\draw[<->]   (4.1,-1.4)-- (6,-1.4) ;
 \draw (-4,-1.8) node[]{$B_{z_1}$};
  \draw (0,-1.8) node[]{$B_{z_2}$};
   \draw (3.2,-1.8) node[]{$B_{z_3}$};
  \draw (5,-1.8) node[]{$B_{z_4}$};
\end{tikzpicture}

\caption{Schematic representation  in dimension $2$ of a function $f$ satisfying the assumptions \textbf{[H1]}, \textbf{[H2]}, and \textbf{[H3]}, and of its restriction $f|_{\pa \Omega}$ to  $\pa \Omega$.  On the figure, $n=4$ and $n_0=2$. }
 \label{fig:H1H2H3H4}

\end{center}
\end{figure}
Under the assumption \textbf{[H2]}, one denotes by~$n_0\in \{1,\ldots,n\}$   the number of global minima of $f|_{\pa \Omega}$, i.e.:
$$f(z_1)= \ldots=f(z_{n_0})<f(z_{n_{n_0+1}})\le\ldots \le f(z_n).$$
On Figure~\ref{fig:H1H2H3H4}, one gives a schematic representation in dimension $2$ of a function~$f$ satisfying the assumptions~\textbf{[H1]},~\textbf{[H2]}, and~\textbf{[H3]}, and of its restriction to~$\pa \Omega$, in the case  $n=4$ and $n_0=2$.
\begin{remark}\label{re.gene-saddle}
The assumption \textbf{[H1]} implies that~$f$ does not have any saddle point (i.e. critical point of index $1$) on $\pa \Omega$. Actually, under \textbf{[H1]}, \textbf{[H2]}, and \textbf{[H3]}, the points~$(z_i)_{i=1,\ldots,n}$ play geometrically the role of saddle points and  are called {\em generalized saddle points} of $f$ on $\pa \Omega$, see~\cite[Section 5.2]{helffer-nier-06}. This  can be explained by the fact that, under  \textbf{[H1]}, \textbf{[H2]}, \textbf{[H3]} and when~$f$ is extended by~$-\infty$ outside~$\overline \Omega$, the points $(z_i)_{i=1,\ldots,n}$ are geometrically   saddle points of $f$ (the extension of~$f$ by $-\infty$ is consistent with the Dirichlet boundary conditions used to define $L^{D,(0)}_{f,h}$) in the following sense. 
For all $i\in \{1,\ldots,n\}$,~$z_i$ is a local minimum of $f|_{\pa \Omega}$ and a  local maximum of~$f|_{D_i}$, where $D_i$ is  the straight line passing through~$z_i$ and orthogonal to  $\pa \Omega$   at~$z_i$. 
\end{remark}

\begin{remark}\label{re.gene-saddle2}
Notice that under  \textbf{[H1]}, \textbf{[H2]}, and \textbf{[H3]}, extending $f$ by reflection outside $\overline \Omega$ in a neighborhood of   $z_i$ also implies that~$z_i$ is  a geometric  saddle point of $f$ as defined in Remark~\ref{re.gene-saddle}. 
 In dimension one, such  a construction was  considered by Kramers in~\cite{Kra} to derive formulas for transition rates, as explained in~\cite{MATK1982}.  Moreover,  as in Remark~\ref{re.facteur2}, it can   be proved in dimension $1$ (exactly as in~\cite[Section A.1.2.2]{BN2017}),  that when reaching $ \pa \Omega\cap \pa \Omega_i $, the probability that the process~\eqref{eq.langevin}    goes in  $\Omega_i$ is  $\frac 12+O(h)$ when  $h\to 0$. To extend this result to  higher dimensions,  one can  use a suitable set of coordinates around~$z_i$. 
\end{remark}
\noindent
Let us now define $g:\overline \Omega\to \mathbb R^+$ by 
\begin{equation}\label{eq.gg}
g(x)=\big \vert \nabla f(x)\big \vert  \,  \text{ when } x\in \Omega \ \text{ and } \ g(x)=\big \vert  \nabla_T f(x) \big \vert  \, \text{ when } x\in \pa\Omega,
\end{equation}
where~$\nabla_Tf$ is the tangential gradient of $f$ in~$\pa \Omega$. Let us recall that for $x\in \pa \Omega$, $\nabla_Tf(x)$ is defined by $\nabla_Tf(x)=\nabla f(x)-(\nabla f(x)\cdot n) \, n$, where $n$ is the unit outward normal  to $\pa \Omega$ at $x$.    The assumptions one needs to state the results in this section depend on the   Agmon distance in~$\overline \Omega$ between the points $(z_i)_{i=1,\ldots,n}$. The  Agmon distance is defined as follows: for any~$x\in\overline \Omega$ and $y\in \overline\Omega$, 
\begin{equation}\label{eq.dag}
d_a(x,y):=\inf_{\gamma \in {\rm Lip}(x,y)}L(\gamma, (0,1)),
\end{equation}
where ${\rm Lip}(x,y)$ is the set of Lipschitz curves~$\gamma: [0,1]\to \overline \Omega$  which are such that $\gamma(0)=x$ and~$\gamma(1)=y$,  and where for $\gamma \in  {\rm Lip}(x,y)$,
$$L(\gamma, (0,1))=\int_0^1 g(\gamma(t))\vert \gamma'(t)\vert dt.$$ 
\begin{remark}
Let us give some common points and differences between  the quasipotential $V$ introduced in~\cite[Section 2]{freidlin-wentzell-84} and  the Agmon distance~\eqref{eq.dag}. Contrary to the quasipotential $V$, the Agmon distance~\eqref{eq.dag} is symmetric. Moreover, let us consider  $x\neq y\in \overline \Omega$  such that there exists a curve  $\gamma:[0,1]\to \overline \Omega$ with $\frac{d}{dt}\gamma(t)=-\nabla f (\gamma(t))$, $\gamma(0) =x$ and  $\gamma(1)=y$. Then, the  Agmon distance~\eqref{eq.dag} between $x$ and $y$ equals $f(x)-f(y)=V(y,x)>0$ but $V(x,y)=0\neq d_a(x,y)$.  
\end{remark}
Finally, let us define the following sets.  For $i\in \{1,\ldots,n\}$, $B_{z_i}$ is the basin of attraction of~$z_i$  for the dynamics $\frac{d}{dt}  x(t)=-\nabla_Tf\big (x(t)\big )$ in~$\pa \Omega$, i.e. $B_{z_i}=\{y\in \pa \Omega,\ \lim_{t\to \infty} x(t)=z_i \text{ if } x(0)=y\}$ (see for instance~Figure~\ref{fig:H1H2H3H4}). 
 Moreover, one defines for $i\in \{1,\ldots,n\}$:
$$B_{z_i}^c:=\pa \Omega\setminus B_{z_i}.$$ 
 
 \subsubsection{Main results.}
 Let us now give the main results of this section.  
 
\begin{proposition}\label{pr1} Let $u_h$  be the eigenfunction associated with the principal eigenvalue $\lambda_h$ of $L^{D,(0)}_{f,h}$ which satisfies normalization~\eqref{norma}. 
Let us assume that the hypotheses \textbf{[H1]}, \textbf{[H2]}, \textbf{[H3]} are satisfied. Then, in the limit $h\to 0$, one has:
\begin{equation}\label{eq.res-lh}
\lambda_h= \frac{\sqrt{ {\rm det \ Hess } f   (x_0) }  }{\sqrt{\pi h}}\sum_{i=1}^{n_0}\frac{  \partial_nf(z_i)    }{    \sqrt{ {\rm det \ Hess } f|_{ \partial \Omega}   (z_i) }  } \ e^{-\frac{2}{h}(f(z_1)-f(x_0))}\left(1+ O(h) \right)
\end{equation}
and
 \begin{equation}\label{eq.res-int-uh}
\int_{\Omega} u_h(x) \ e^{- \frac{2}{h} f(x) }dx= \frac{ \pi ^{\frac{d}{4} } }{  \left({\rm det \ Hess } f   (x_0)  \right)^{1/4}  } \  h^{\frac d4} \ e^{-\frac{1}{h}f(x_0)} (1+O(h) ).
 \end{equation}
\end{proposition}
Furthermore, one obtains  the following theorem on the asymptotic behaviour of~$\pa_nu_h$, which is one of the main results of~\cite{DLLN}. 

\begin{theorem}\label{th1}
Let us assume that \textbf{[H1]}, \textbf{[H2]}, and \textbf{[H3]} are satisfied  and that the following inequalities hold:
\begin{equation}\label{eq.da1}
f(z_1)-f(x_0)> f(z_n)-f(z_1)
\end{equation}
and for all $i\in \{1,\ldots,n\}$, 
\begin{equation}\label{da2}
d_a(z_i,B_{z_i}^c)>\max[f(z_n)-f(z_i),\, f(z_i)-f(z_1)].
\end{equation}
Let~$i\in \{1,\ldots,n\}$ and~$\Sigma_i\subset \pa \Omega$  be an open set containing $z_i$ and such that $\overline{\Sigma_i} \subset B_{z_i}$. Let $u_h$  be the eigenfunction associated with the principal eigenvalue of~$L^{D,(0)}_{f,h}$ which satisfies~\eqref{norma}.  Then, in the limit~$h\to 0$:
\begin{equation}\label{pa-nuh}
\int_{\Sigma_i}   \partial_n u_h \,  e^{- \frac{2}{h}  f}  =C_i(h) \,e^{-\frac{2f(z_i)-f(x_0)}{h}}\left(1+ O(h) \right),
\end{equation}
where~$C_i(h)=-\frac{\displaystyle  ({\rm det \ Hess } f   (x_0))^{1/4}   \partial_nf(z_i)   2\pi^{\frac{d-2}{4}}    }{    \displaystyle\sqrt{ {\rm det \ Hess } f|_{\partial \Omega}   (z_i) }  }\, h^{\frac{d-6}{4}}$. 
\end{theorem}
These results have the following consequences. 
\begin{corollary}\label{cc1}
Let us assume that all the assumptions of Theorem~\ref{th1} are satisfied. Let~$i\in \{1,\ldots,n\}$ and~$\Sigma_i\subset \pa \Omega$ be  an open set containing $z_i$ and such that $\overline{\Sigma_i} \subset B_{z_i}$. Then, in the limit $h\to 0$:
\begin{align}\nonumber
\mathbb P_{\nu_h} \left[ X_{\tau_{\Omega}} \in \Sigma_i\right]&=  \frac{\partial_nf(z_i)  }{ \sqrt{ {\rm det \ Hess } f|_{\partial \Omega }   (z_i) }}  
\left (\sum_{k=1}^{n_0} \frac{ \partial_nf(z_k) }{\sqrt{ {\rm det \ Hess } f|_{\partial \Omega }   (z_k) } } \right)^{-1}\\
\label{Proba-nuh}
&\quad \times  e^{-\frac{2}{h} (f(z_i)-f(z_1))}  (  1+   
 O(h)  ),
 \end{align} 
where $\nu_h$ is the quasi-stationary distribution of the process~\eqref{eq.langevin} in~$\Omega$ (see~\eqref{exp.qsd}). Moreover,  if $\Sigma_i$ is the common boundary between the state $\Omega$ and a state~$\Omega_i$, then, when~$h\to 0$ 
\begin{equation}
\label{eq.kii}
k_{i}^L=\frac{1}{\sqrt{\pi h}}  \partial_nf(z_i)  \frac{ \sqrt{  {\rm det \ Hess } f  (x_0) } }{  \sqrt{  {\rm det \ Hess } f|_{\partial \Omega }   (z_i) }} e^{-\frac{2}{h}(f(z_i)-f(x_0))} (1+ O (h)),
\end{equation}
where $k_{i}^L$ is the transition rate~\eqref{kijl} to go from $\Omega$ to $\Omega_i$. 
\end{corollary}
\begin{sloppypar}
Notice that since $z_i$ is not a  saddle  point of~$f$, the prefactor in~\eqref{eq.kii} is not the prefactor $\frac 12 A_i$ (see Remark~\ref{re.gene-saddle2} for the explanation of the multiplicative term~$\frac 12$), where $A_i$ is defined by~\eqref{eq.prefactor}, but it is  actually  the expected prefactor  for a generalized saddle point of $f$ (see Remarks~\ref{re.gene-saddle} and~\ref{re.gene-saddle2}).

The asymptotic estimate~\eqref{Proba-nuh}  is a consequence of Proposition~\ref{pr1}, Theorem~\ref{th1} together with~\eqref{expression-P}, and~\eqref{eq.kii} is a consequence of Proposition~\ref{pr1}, Theorem~\ref{th1} and~\eqref{expression-kijl}.
 The main difficulty is to prove~\eqref{pa-nuh} which requires a sharp equivalent of the quantity~$\int_{\Sigma_i}   \partial_n u_h \,  e^{- \frac{2}{h}  f} $ when $z_i$ is not a global minimum of $f$ on $\pa \Omega$, i.e. when~$i\in\{n_0+1,\ldots,n\}$.  \\
In~\cite{DLLN}, numerical simulations are provided to check that~\eqref{Proba-nuh} holds and to discuss the necessity of the assumptions~\eqref{da2} to obtain~\eqref{Proba-nuh}. Furthermore, in~\cite{DLLN},  the results~\eqref{pa-nuh} and~\eqref{Proba-nuh} are generalized to sets~$\Sigma\subset \pa \Omega$ which do not necessarily contain a point~$z\in \{z_1,\ldots,z_n\}$: this is the other main results of~\cite{DLLN} which is not presented here.  Moreover, with the help of  ``leveling" results on the function~$x\mapsto \mathbb E_x[F(X_{\tau_\Omega})]$, we generalized~\eqref{Proba-nuh}  to deterministic initial conditions in~$\Omega$ (i.e. when $X_0=x\in \Omega$) which are the initial conditions considered in the theory of large deviations~\cite{freidlin-wentzell-84}. 
\end{sloppypar}
The proofs of Proposition~\ref{pr1} and Theorem~\ref{th1} are based on tools from semi-classical analysis and more precisely, they are based on techniques developed in~\cite{HeSj4,helffer-sjostrand-85b, helffer-sjostrand-85a,helffer-sjostrand-84, helffer-nier-06, le-peutrec-10}.

 \subsubsection{Starting points of the proofs of Proposition~\ref{pr1} and Theorem~\ref{th1}.}
Let  us recall that $u_h$ is the eigenfunction associated with  the  principal eigenvalue~$\lambda_h$ of $L^{D,(0)}_{f,h}$ which satisfies normalization~\eqref{norma}. In view of~\eqref{expression-P} and in order to obtain~\eqref{Proba-nuh}, one wants to study the asymptotic behaviour when $h\to 0$ of~$\nabla u_h$  on~$\pa\Omega$.
The starting point of the proofs of Proposition~\ref{pr1} and Theorem~\ref{th1} is the fact that  $\nabla u_h$ is solution to an eigenvalue problem   for the same eigenvalue~$\lambda_h$. Indeed, recall that $u_h$ is  solution to $
 L^{(0)}_{f,h}\, u_h  =  \lambda_h u_h$ in $\Omega$ and  $u_h = 0 $ on $\partial \Omega$. If one differentiates this relation, $\nabla u_h$ is solution to 
\begin{equation}\label{eq:L1-intro}
\left\{
\begin{aligned}
L^{(1)}_{f,h} \nabla u_h &= \lambda_h \nabla u_h \text{ in $\Omega$,}\\
\nabla_T u_h& = 0 \text{ on $\partial \Omega$,}\\
\left(-\frac{h}{2} {\rm div} + \nabla f \cdot \right) \nabla u_h & = 0 \text{ on $\partial \Omega$,}\\
\end{aligned}
\right.
\end{equation}
where  $
L^{(1)}_{f,h}=- \frac{h}{2} \Delta + \nabla f \cdot \nabla + {\rm
  Hess} \, f$ 
is an operator acting on $1$-forms (namely on vector fields). In the following the operator $L^{(1)}_{f,h}$ with tangential boundary conditions~\eqref{eq:L1-intro} is denoted by~$L^{D,(1)}_{f,h}$. From~\eqref{eq:L1-intro},  $\nabla u_h$ is therefore an eigenform of $L^{D,(1)}_{f,h}$ associated with $\lambda_h$. 
For $p\in\{0,1\}$, let us denote, by  $ \pi_h^{(p)}$  the orthogonal projector  of $L^{D,(p)}_{f,h}$ associated with the eigenvalues of $L^{D,(p)}_{f,h}$ smaller than $\frac{\sqrt h}{2}$. Another crucial ingredient for the proofs of Proposition~\ref{pr1} and Theorem~\ref{th1} is the fact that, from~\cite[Chapter 3]{helffer-nier-06}, 
\begin{equation}\label{eq.dim-no}
\Ran \, \pi_h^{(0)} =  {\rm Span}\, u_h  \text{ and}  \   {\rm dim} \,  \Ran\, \pi_h^{(1)} =n.
\end{equation} 
Therefore, from~\eqref{eq:L1-intro}, it holds 
\begin{equation}\label{eq.nu-in}
\nabla u_h\in \Ran\, \pi_h^{(1)},  
\end{equation}
and from~\eqref{norma} and  the fact that $\langle L^{(0)}_{f,h}\, u_h, u_h\rangle_{L^2_w}=\frac h2\Vert \nabla u_h \Vert_{L^2_w}^2$, one has
\begin{equation}\label{eq.nu-in-l}
\lambda_h=\frac h2\Vert \nabla u_h \Vert_{L^2_w}^2.
\end{equation}
Thus, to study the asymptotic behaviour when $h\to 0$ of  $\lambda_h$, $u_h$ and $\nabla u_h$, we construct a suitable  orthonormal basis of  $\Ran\, \pi_h^{(1)} $. This basis is constructed using so-called quasi-modes. 
\subsubsection{Sketch of the proofs of Proposition~\ref{pr1} and Theorem~\ref{th1}.}

Let us give the sketch of the proof of~\eqref{Proba-nuh} which is the main result of~\cite{DLLN}. Recall that from Proposition~\ref{fried-bis},  one works in the Hilbert space $L^2_w(\Omega)$. The spaces $L^2_w(\Omega)$ and $H^{1}_w(\Omega)$   (see~\eqref{eq.H1w})  extend naturally on $1$-forms  as follows 
$$\Lambda^1L^2_w(\Omega):=\Big \{u=\, ^t(u_1,\ldots,u_d):\Omega \to \mathbb R^d,  \, \forall k\in\{1,\ldots,d\},\  \int_\Omega u_k^2
 e^{-\frac{2}{h} f} < \infty\Big\},$$
 and 
$$\Lambda^1 H^{1}_w(\Omega):= \big\{u=\, ^t(u_1,\ldots,u_d):\Omega \to \mathbb R^d,  \, \forall (i,k)\in\{1,\ldots,d\}^2,\    \pa_i u_k \in L^2_w(\Omega)\big  \}.$$
In the following, one denotes by  $\Vert.\Vert_{L^2_w}$ (resp. $\Vert.\Vert_{H^{1}_w}$) the norm of  $L^2_w(\Omega)$ and  of $\Lambda^1L^2_w(\Omega)$ (resp. $H^1_w(\Omega)$ and~$\Lambda^1 H^{1}_w(\Omega)$). Finally,~$\langle.,.\rangle_{L^2_w}$ stands for both  the scalar product associated with the norm of~$L^2_w(\Omega)$  and with the norm  of $\Lambda^1L^2_w(\Omega)$.
In view of~\eqref{eq.nu-in} and~\eqref{eq.dim-no}, one has for all orthonormal basis  $(\psi_j)_{j\in \{1,\ldots,n\}}$ of~$\Ran\, \pi_h^{(1)}$, in $L^2_w(\Omega)$:
\begin{equation}\label{eq.dec-pauh1}
\nabla u_h =\sum_{j=1}^n \langle \nabla u_h, \psi_j \rangle_{L^2_w}  \psi_j,
\end{equation}
and from~\eqref{eq.nu-in-l}, it holds
\begin{equation}\label{eq.lambda-dec}
\lambda_h =\frac h2 \sum_{j=1}^n \big \vert \langle \nabla u_h, \psi_j \rangle_{L^2_w}  \big \vert^2.
\end{equation}
In particular, one has for all $k\in  \{1,\ldots,n\}$,
\begin{equation}\label{eq.dec-pa-uh}
\int_{\Sigma_k}\pa_n u_h \, e^{-\frac 2h f} =\sum_{j=1}^n \langle \nabla u_h, \psi_j \rangle_{L^2_w}    \int_{\Sigma_k}   \psi_j\cdot n\,  e^{-\frac 2h f},
\end{equation}
where we recall that $\Sigma_k$ is  an open set of~$\pa \Omega$ such that  $z_k\in \Sigma_k$ and $ \overline{\Sigma_k}\subset B_{z_k}$. 
\medskip

\noindent
\textit{Step 1: approximation of $u_h$}. Under  \textbf{[H1]}, \textbf{[H2]}, and \textbf{[H3]}, it is not difficult to find a good  approximation of $u_h$.  Indeed, let us consider, 
\begin{equation}\label{tilde-u}
\tilde u:= \frac{\chi }{\Vert \chi\Vert_{L^2_w}},
\end{equation}
where $\chi\in C^\infty_c(\Omega,\mathbb R^+)$ and $\chi=1$ on $\{x\in \Omega, d(x,\pa \Omega )\ge \varepsilon\}$ where $\varepsilon >0$. In particular,  for $\varepsilon$ small enough, $\chi=1$ in a neighboorhood of  
$x_0$ (which is assumed in the following). Let us explain why $\tilde u$ is a good approximation of~$u_h$. Since $L^{D, (0)}_{f,h}$ is self adjoint on  $L^2_w(\Omega)$, one has  
$$\big \Vert (1-\pi_h^{(0)})\tilde u \big \Vert_{L^2_w}^2 \le \frac{C}{\sqrt h} \big \langle L^{D, (0)}_{f,h}\tilde u, \tilde u \big\rangle_{L^2_w}= \frac{Ch}{2\sqrt h}\frac{  {\int_{\Omega} \vert \nabla \chi\vert^2 e^{-\frac 2h f}}}{  {\int_{\Omega} \chi ^2 e^{-\frac2h f}} }.$$
Since $f(x_0)=\min_{\Omega}f<\min_{\pa \Omega}f$ and $x_0$  is the unique global minimum of $f$ on $\overline \Omega$ (see~\textbf{[H2]}), one has using Laplace's method ($x_0$ is a non degenerate critical point of  $f$ and $\chi(x_0)=1$):
 $$\int_{\Omega} \chi ^2 e^{-\frac2h f}=\frac{   (\pi h)^{\frac d2}   }{ \sqrt{{\rm det} {\rm Hess} f(x_0)}    } e^{-\frac 2h f(x_0)} (1+O(h)).$$
Therefore, for any $\delta>0$, choosing  $\varepsilon$ small enough, it holds when  $h\to 0$: $$\big \Vert (1-\pi_h^{(0)})\tilde u \big \Vert_{L^2_w}^2=O(e^{-\frac 2h (f(z_1)-f(x_0)-\delta)}),$$
and thus:
$$\pi_h^{(0)}\tilde u=\tilde u + O(e^{-\frac 1h (f(z_1)-f(x_0)-\delta)}) \, \text{ in } \,  L^2_w(\Omega).$$ 
 From~\eqref{eq.dim-no} and since $\chi\ge 0$,   one has for any $\delta>0$ (choosing  $\varepsilon$ small enough),  when $h\to 0$ 
\begin{equation}\label{=uh=}
u_h=\frac{  \pi_h^{(0)} \tilde u } { \Vert  \pi_h^{(0)} \tilde u \Vert_{L^2_w}}=  \tilde u \,   +O(e^{-\frac 1h (f(z_1)-f(x_0)-\delta)}) \, \text{ in } \,  L^2_w(\Omega). 
\end{equation}
Since $ \Vert  \tilde u  \Vert_{L^2_w}=1$, this last relation  justifies that  $\tilde u$  is a good approximation of $u_h$ in $L^2_w(\Omega)$. Notice that~\eqref{=uh=}   implies~\eqref{eq.res-int-uh}.  
  \medskip
 
\noindent
\textit{Step 2: construction  of a basis of $\Ran\, \pi_h^{(1)}$ to prove Theorem~\ref{th1}}. In view of~\eqref{eq.dec-pa-uh}, the idea is to construct a family of $1$-forms  $(\widetilde \psi_j)_{j\in \{1,\ldots,n\}}$ which forms, when projected on  $\Ran \,  \pi_h^{(1)}$,   a basis of $\Ran \,  \pi_h^{(1)}$  and which allows to obtain sharp asymptotic estimates on~$\partial_n u_h$ on all the $\Sigma_j$'s when $h\to 0$.  In the literature, such a $1$-form $\widetilde \psi_j$   is called a quasi-mode (for~$L^{D,(1)}_{f,h}$).  A quasi-mode for~$L^{D,(1)}_{f,h}$  is a smooth $1$-form~$w$ such that for some norm, it holds when $h\to 0$:
\begin{equation}\label{eq.residue}
 \pi_h^{(1)}w=w +o(1), 
 \end{equation}
 \begin{sloppypar}
 \noindent
To prove Theorem~\ref{th1}, one of the major issues   is the construction of a basis~$(\widetilde \psi_j)_{j\in \{1,\ldots,n\}}$ so that the remainder term $o(1)$ in~\eqref{eq.residue}, when~$w=\widetilde \psi_k$, is of the order (see~\eqref{da2})
\begin{equation}\label{eq.pi-i}
\big \Vert (1-\pi_h^{(1)})\widetilde \psi_k \big \Vert_{H^1_w}=O \big (e^{-\frac 1h \max[f(z_n)-f(z_k),\, f(z_k)-f(z_1)]}\big).
\end{equation}
\end{sloppypar}
\noindent
 This implies  that $
\big (\pi_h^{(1)}\widetilde \psi_j \big)_{j\in \{1,\ldots,n\}} 
$ is a basis of $\Ran \,  \pi_h^{(1)}$ and  above all,  after   a Gram-Schmidt procedure on~$\big (\pi_h^{(1)}\widetilde \psi_j \big)_{j\in \{1,\ldots,n\}}$,   when $h\to 0$,  that for all $k\in  \{1,\ldots,n\}$ (see~\eqref{eq.dec-pa-uh}):
 \begin{equation}\label{eq.dec-pa-uh-bis}
  \int_{\Sigma_k}\pa_n u_h \, e^{-\frac 2h f}= \sum_{j=1}^n  \langle \nabla \tilde u, \widetilde \psi_j \rangle_{L^2_w}    \int_{\Sigma_k}  \widetilde  \psi_j\cdot n\,  e^{-\frac 2h f} +O\big ( e^{-\frac{2f(z_k)-f(x_0) +c }{h}} \big)
\end{equation}
and (see~\eqref{eq.lambda-dec})
 \begin{equation}\label{eq.lambda-dec-bis}
\lambda_h =\frac h2 \sum_{j=1}^n \vert \langle \nabla \tilde u, \widetilde \psi_j \rangle_{L^2_w}  \vert^2 + O\big ( e^{-\frac{2}{h}(f(z_1)-f(x_0)+c) } \big)
\end{equation}
for some  $c>0$   independent of $h$. Here,  we recall, $\tilde u$ (see~\eqref{tilde-u}) is a   good approximation of $u_h$ (see~\eqref{=uh=}).
Let us now explain  how we will construct the family 
$\big (\widetilde \psi_j \big)_{j\in \{1,\ldots,n\}}$in order to obtain~\eqref{eq.dec-pa-uh-bis} and~\eqref{eq.lambda-dec-bis}. Then, we explain how the terms $\Big(\displaystyle \int_{\Sigma_j}  \widetilde  \psi_j\cdot n\,  e^{-\frac 2h f}\Big)_{j\in \{1,\ldots,n\}}$ and $\Big(\langle \nabla \tilde u, \widetilde \psi_j \rangle_{L^2_w}  \Big)_{j\in \{1,\ldots,n\}}$ appearing in~\eqref{eq.dec-pa-uh-bis} and~\eqref{eq.lambda-dec-bis}  are computed.
\medskip

\begin{sloppypar}
\noindent
\textit{Step 2a: construction of the family $(\widetilde \psi_j)_{j\in \{1,\ldots,n\}}$}. 
 To construct each $1$-form~$\widetilde \psi_j$, the idea is to construct an operator $L^{(1)}_{f,h}$  with mixed tangential Dirichlet and Neumann boundary conditions on a domain $\dot{\Omega}_j\subset \Omega$ which is such that $\big  (\{z_1,\ldots,z_n\} \cup \{x_0\} \big )\cap \dot{\Omega}_j=\{z_j\}$.   For $j\in\{1,\ldots,n\}$, $\widetilde  \psi_j$ is said to be associated with the generalized saddle point $z_j$.
  The goal of the boundary conditions is to ensure that when $h\to 0$, each of these  operators  has only  one exponentially small eigenvalue (i.e. this eigenvalue is  $O\big(e^{-\frac ch}\big)$ for some $c>0$ independent of $h$), the other eigenvalues being larger than $\sqrt h$. Then, we show that each of these small eigenvalues actually equals~$0$ using the Witten complex structure associated  with these boundary conditions on $\pa \dot \Omega_j$. To construct such operators $L^{(1)}_{f,h}$ with mixed boundary conditions on~$\dot{\Omega}_j$, the recent results of~\cite{jakab-mitrea-mitrea-09} and~\cite{goldshtein-mitrea-mitrea-11} are used.
The~$1$-form~$\widetilde \psi_j$ associated with $z_j$ is then defined using an eigenform~$v_{h,j}^{(1)}$ associated with the eigenvalue $0$ of the operator $L^{(1)}_{f,h}$   associated with mixed boundary conditions on~$\dot{\Omega}_j$:  
 \begin{equation}\label{tilde-psij}
 \widetilde \psi_j:= \frac{ \chi_j \,  v_{h,j}^{(1)}    }{   \Vert \chi_j  \, v_{h,j}^{(1)} \Vert_{L^2_w}  },
\end{equation}
where $\chi_j$ is a well chosen cut-off function  with support in $\overline{\dot \Omega_j}$.  Notice that for~$j\in \{1,\ldots,n\}$, the quasi-mode $\widetilde \psi_j$ is not only constructed in a neighbourhood of $z_j$: it has a support  as large as needed in $\Omega$. This is a difference with previous construction in the literature, such as~\cite{helffer-nier-06}.      We need such quasi-modes for the following reasons. Firstly, we   compute the probability that the process~\eqref{eq.langevin} leaves $\Omega$ through open sets $\Sigma_j$ which are arbitrarily large in $B_{z_j}$. Secondly, we use the fact that the quasi-mode~$\widetilde \psi_j$ decreases very fast away from~$z_j$ to get~\eqref{eq.pi-i}. This is needed to state the hypothesis~\eqref{da2} in terms of Agmon distances, see next step. 
 \end{sloppypar}
\medskip

\noindent
\textit{Step 2b: Accuracy of the quasi-mode $\widetilde \psi_j$ for ${j\in \{1,\ldots,n\}}$}. 
To obtain a sufficiently small remainder term in~\eqref{eq.residue} (to get~\eqref{eq.pi-i} and then~\eqref{eq.dec-pa-uh-bis}), one needs to quantify the decrease of the quasi-mode $\widetilde \psi_j$ outside a neighboorhood of $z_j$. This decrease   
is obtained with Agmon estimates on $v_{h,j}^{(1)}$ which allow to localize~$\widetilde \psi_j$ in a neighboorhood of $z_j$. 
For ${j\in \{1,\ldots,n\}}$, we prove the following Agmon estimate on $v_{h,j}^{(1)}$:
\begin{equation}\label{eq.est-ag}
\big \Vert \chi_j \,  v_{h,j}^{(1)} e^{\frac 1h d_a(.,z_j) }\big  \Vert_{H^1_w}=O(h^{-N}),
\end{equation}
for some $N\in \mathbb N$ and where $d_a$ is the Agmon distance defined in~\eqref{eq.dag}.  To obtain~\eqref{eq.est-ag}, we study the properties of this distance. The boundary of~$\Omega$   introduces technical difficulties. The Agmon estimate~\eqref{eq.est-ag} is obtained adapting to our case techniques developed in~\cite{helffer-nier-06, le-peutrec-10}.
For all~$j\in \{1,\ldots,n\}$, using the fact that~$\big \Vert (1-\pi_h^{(1)})\widetilde \psi_j \big \Vert_{L^2_w}^2 \le \frac{C}{\sqrt h} \big \langle L^{D, (1)}_{f,h}\widetilde \psi_j, \widetilde \psi_j\big\rangle_{L^2_w}$ and~\eqref{eq.est-ag}, 
one shows that  
$$\big \Vert (1-\pi_h^{(1)})\widetilde \psi_j \big \Vert_{L^2_w}^2 \le C\, h^{-q} \, e^{-\frac 2h \inf_{\text{supp}\nabla \chi_j} d_a(.,z_j)  },$$
for some $q>0$. 
Thus, in order to get~\eqref{eq.pi-i}, the support of $\nabla \chi_j$ has to be arbitrarily close to~$x_0$ and~$B_{z_j}^c$. This explains the assumptions~\eqref{eq.da1} and~\eqref{da2}, and the fact that the quasi-mode $\widetilde \psi_j$ is not constructed in a neighboorhood of $z_j$ but in a domain~$\dot \Omega_j$ arbitrarily large in~$\Omega$. This is one of the main differences compared with~\cite{helffer-nier-06}.  
At the end of this step, one has a family~$(\widetilde \psi_j)_{j\in \{1,\ldots,n\}}$ which satisfies~\eqref{eq.pi-i}. This allows us to obtain, in the limit $h\to 0$   (see~\eqref{eq.dec-pa-uh-bis}), for some $c>0$ independent of~$h$ and for all $  k\in  \{1,\ldots,n\}$:
$$
  \int_{\Sigma_k}\pa_n u_h \, e^{-\frac 2h f}=   \sum_{j=1}^n\langle \nabla \tilde u, \widetilde \psi_j \rangle_{L^2_w}    \int_{\Sigma_k}  \widetilde  \psi_j\cdot n\,  e^{-\frac 2h f} +O\big ( e^{-\frac{2f(z_k)-f(x_0)+c}{h}} \big).
$$
 
 \medskip
 \noindent
\textit{Etape 3: computations of $\Big(\displaystyle \int_{\Sigma_j}  \widetilde  \psi_j\cdot n\,  e^{-\frac 2h f}\Big)_{j\in \{1,\ldots,n\}}$ and $\Big(\langle \nabla \tilde u, \widetilde \psi_j \rangle_{L^2_w}  \Big)_{j\in \{1,\ldots,n\}}$}.
In view of~\eqref{eq.dec-pa-uh-bis} and~\eqref{eq.lambda-dec-bis}, for all $j\in \{1,\ldots,n\}$,  one needs to compute the terms $$  \int_{\Sigma_j}   \widetilde \psi_j\cdot n\,  e^{-\frac 2h f} \text{ and } \langle \nabla \tilde u, \widetilde \psi_j \rangle_{L^2_w}.$$
To do that, we use for all  ${j\in \{1,\ldots,n\}}$ a WKB approximation of $v_{h,j}^{(1)}$, denoted by $v_{z_j,wkb}^{(1)}$.   In the literature we follow, $v_{z_j,wkb}^{(1)}$  is constructed in a neighboorhood of $z_j$  (see~\cite{helffer-nier-06, le-peutrec-10}). To prove Theorem~\ref{th1}, we extend the construction of~$v_{z_j,wkb}^{(1)}$ to  neighbourhoods in $\overline \Omega$  of arbitrarily large closed   sets included in~$B_{z_j}$ (indeed, there is no restriction on the size of $\Sigma_j$ in $B_{z_j}$).  
Then, the comparison between~$v_{h,j}^{(1)}$ and~$v_{z_j,wkb}^{(1)}$  is also extended to  neighbourhoods in $\overline \Omega$  of arbitrarily large closed   sets included in $B_{z_j}$.  
Once the terms $\Big(\int_{\Sigma_j}  \widetilde  \psi_j\cdot n\,  e^{-\frac 2h f}\Big)_{j\in \{1,\ldots,n\}}$ and $\Big(\langle \nabla \tilde u, \widetilde \psi_j \rangle_{L^2_w}  \Big)_{j\in \{1,\ldots,n\}}$ are computed, one concludes the proof of~\eqref{eq.res-lh} using~\eqref{eq.lambda-dec-bis} and the proof of~\eqref{pa-nuh}  using~\eqref{eq.dec-pa-uh-bis}. 



\subsection{Most probable   exit points from a bounded domain} 
\subsubsection{Setting and motivation.}
In this section, we present recent results from~\cite{di-gesu-lelievre-le-peutrec-nectoux-16} on the concentration of the law of $X_{\tau_\Omega}$ on  a subset of ${\rm argmin}_{\pa \Omega}f=\{z\in \pa \Omega, \ f(z)=\min_{\pa \Omega}f\}$ when $h\to 0$ in a more general geometric setting than the one of Theorem~\ref{th1}. The main purpose of these results is to prove an asymptotic formula when $h\to 0$ for the concentration of the law of $X_{\tau_\Omega}$ on a set of points of ${\rm argmin}_{\pa \Omega} f$ when~$\Omega$ contains several local minima of $f$ and when $\pa_n f$ is not necessarily positive on~$\pa \Omega$. 

Let $\mathcal Y \subset \pa \Omega$.  We say that the law of~$X_{\tau_{\Omega}}$ concentrates on~$\mathcal Y$ if for all neighborhood~$\mathcal  V_{\mathcal Y}$ of~${\mathcal Y}$ in $\pa \Omega$, one has 
~$$\lim \limits_{h\to 0}\mathbb P  \left [ X_{\tau_{\Omega}} \in \mathcal  V_{\mathcal Y}\right]=1,$$ 
and if for all $x \in \mathcal Y$ and all  neighborhood~$\mathcal  V_x$ of~$x$ in~$\pa \Omega$ , it holds:  
~$$\lim \limits_{h\to 0}\mathbb P  \left [ X_{\tau_{\Omega}} \in \mathcal  V_x\right]>0.$$ 
In~\cite{MS77,schuss90, MaSc}, when~$\pa_n f(x)=0$ for all~$x\in \pa \Omega$ or when~$\pa_n f(x)>0$ for all~$x\in \pa \Omega$ (and with additional assumptions on $f$), it has been shown that the law of~$X_{\tau_{\Omega}}$ concentrates on points where~$f$ attains its minimum on~$\pa \Omega$ (see~\eqref{eq.MS}).  Later on, it has been proved in~\cite{kamin1979elliptic,kamin-78,perthame-90, day1984a,day1987r} when~$\pa_nf>0$ on~$\pa \Omega$ and~$f$  has a unique non degenerate critical point in~$\Omega$ (which is necessarily its global minimum in~$\overline \Omega$). Tools developed in semi-classical analysis allow us to generalize this geometric setting. For instance, we consider  several critical points of~$f$ in~$\Omega$  and we drop  the assumptions~$\pa_nf>0$ on~$\pa \Omega$ (however we do not consider the case when~$f$  has saddle points on~$\partial
\Omega$). Assuming that~$f$ and~$f|_{\pa \Omega}$ are Morse functions, and~$\vert \nabla f \vert \neq 0$ on~$\pa \Omega$,  we raise the following questions:
\begin{itemize}
\item  What are the geometric conditions ensuring that, when~$X_0\sim \nu_h$,   the law  of~$X_{\tau_{\Omega}}$ concentrates on  points where~$f$ attains its minimum on~$\pa \Omega$ 
(or a subset of these points)?

\item \begin{sloppypar}
What are the conditions which ensure that these results extend to some deterministic initial conditions in 
$\Omega$ ?
 \end{sloppypar}
\end{itemize}
The results of~\cite{di-gesu-lelievre-le-peutrec-nectoux-16} aim at answering these questions. Let us recall that when~$f$ and~$f|_{\pa \Omega}$ are Morse functions and when~$\vert \nabla f \vert \neq 0$ on~$\pa \Omega$,  the  elements of the set
\begin{equation}
\label{PS-gene}
\{ z \text { is a local minimim of }  f|_{\pa \Omega}\}\cap\{z\in \pa \Omega, \, \pa_nf(z)>0\}
\end{equation} 
are the generalized saddle points of~$f$ on $\pa \Omega$ and play the role of saddle points of~$f$ on~$\pa \Omega$, see Remark~\ref{re.gene-saddle}. 
 Before stating the main results of~\cite{di-gesu-lelievre-le-peutrec-nectoux-16}, let us discuss the two  questions above
  with one-dimensional examples.
  
\begin{remark}
The assumption that the drift term $b$ in~\eqref{eq.langevin} is of the form $b=-\nabla f$  is essential here to the existence of a limiting exit distribution of $\Omega$ when $h\to 0$. If it is not the case and when for instance  the boundary of $\Omega$  is a periodic orbit of the dynamics $\frac{d}{dt}  x(t)=b\big (x(t)\big )$,  the phenomenon of cycling discovered by Day in~\cite{day1992,day1996} prevents the
existence of a limiting exit distribution when $h\to 0$.    We also refer to~\cite{berglund2014noise,berglund2004noise2,berglund2014noise3} for the study of this phenomenon of cycling. 
\end{remark}

\subsubsection{One-dimensional examples.}
To discuss  the two  questions raised in the previous section, one considers two one-dimensional examples. 

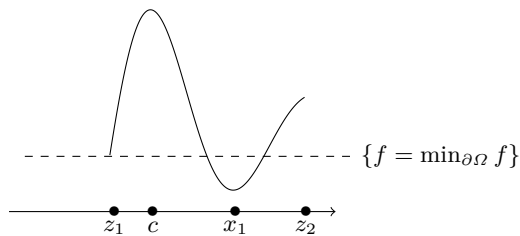
\begin{figure}
\begin{center}
\begin{tikzpicture}[scale=0.7]
  \draw [domain=-0.17:2.3*pi, scale=0.5, samples=300, smooth] plot (\x,{6*sin(\x r)*exp(-0.2*\x)}); 
  \draw[->] (-2,-1.6) -- (4.2,-1.6);
  \draw (2.3,-1.9) node[]{$x_1$};
  \draw (2.3,-1.6) node[]{$\bullet$};
  \draw [dashed] (-1.7,-0.55)--(4.6,-0.55); 
   \draw (6.2,-0.55) node[]{$\{f=\min_{\pa \Omega}f\}$};
   \draw (0,-1.9) node[]{$z_1$};
  \draw (0,-1.6) node[]{$\bullet$};  
  \draw (3.64,-1.9) node[]{$z_2$};
  \draw (3.64,-1.6) node[]{$\bullet$};
      \draw (0.73,-1.9) node[]{$c$};
  \draw (0.73,-1.6) node[]{$\bullet$};
  
    \end{tikzpicture}
    \caption{ Example of a function~$f$ such that, starting  from the global minimum~$x_1$ of~$f$ in~$\Omega$ or from the quasi-stationary distribution~$\nu_h$,  the law of~$X_{\tau_{\Omega}}$  concentrates on $z_2$ whereas  $f(z_2)>\min_{\pa \Omega} f=f(z_1)$ . }
    \label{fig:ex1-1}
    \end{center}
\end{figure}

\medskip
\noindent
\textit{Example 1}. 
The goal  is here to construct a one-dimensional example for which, starting from the global minimum of~$f$ in~$\Omega$ or from the quasi-stationary distribution~$\nu_h$,  the law of~$X_{\tau_{\Omega}}$ does not concentrate on points  where~$f$  attains its minimum on~$\pa \Omega$.  To this end, let us consider the function~$f$ represented in  Figure~\ref{fig:ex1-1} for which one has the following result.

\begin{proposition}\label{pr-ex1}
Let~$z_1<z_2$ and~$f\in C^{\infty}([z_1,z_2], \mathbb R)$    be a Morse function. Let us assume that~$f(z_1)<f(z_2)$,~$\{x\in [z_1,z_2], f'(x)=0\}=\{c,x_1\}$ with~$z_1<c<x_1<z_2$ and~$f(x_1)<f(z_1)<f(z_2)<f(c)$ (see~Figure~\ref{fig:ex1-1}). Then, for all~$x\in (c,z_2]$, there exists~$\varepsilon>0$ such that when~$h\to 0$:
$$
\mathbb P_{x}  [ X_{\tau_{(z_1,z_2)}}=z_1 ]=O (e^{- \frac{\varepsilon}{h} }) \text{ and thus }    \mathbb P_{x}  [ X_{\tau_{(z_1,z_2)}}=z_2 ]= 1+O (e^{- \frac{\varepsilon}{h} }).
$$
Moreover, there exists~$\varepsilon>0$ such that when~$h\to 0$:
$$
\mathbb P_{\nu_h}  [ X_{\tau_{(z_1,z_2)}}=z_1 ]=O (e^{- \frac{\varepsilon}{h} })   \ \text{  and thus }  \ \mathbb P_{\nu_h}  [ X_{\tau_{(z_1,z_2)}}=z_2 ]=1+O (e^{- \frac{\varepsilon}{h} }),
 $$
where~$\nu_h$ is the quasi-stationary distribution of the process~\eqref{eq.langevin} in~$(z_1,z_2)$. 
\end{proposition}
The proof of Proposition~\ref{pr-ex1} is based on the fact that in one dimension, explicit formulas can be written for $x\mapsto \mathbb P_x [ X_{\tau_{(z_1,z_2)}}=z_j ]$ ($j\in \{1,2\}$), see~\cite[Section A.5.3.1]{BN2017} or~\cite{di-gesu-lelievre-le-peutrec-nectoux-16}.  
According to Proposition~\ref{pr-ex1}, when~$h\to0$ and when~$X_0=x\in (c,z_2)$ or~$X_0\sim\nu_h$, the process~\eqref{eq.langevin} leaves~$\Omega=(z_1,z_2)$ through~$z_2$. However, the generalized saddle point~$z_2$ (see~\eqref{PS-gene}) is not the global minimum of~$f$  on~$\pa \Omega$. 
This fact can be explained as follows: the potential  barrier~$f(c)-f(x_1)$ is larger than  the potential  barrier~$f(z_2)-f(x_1)$. 
Thus, the law of~$X_{\tau_\Omega}$ when~$X_0=x\in (c,z_2)$  cannot concentrate on~$z_1$  since it is  less costly to leave~$\Omega$ through~$z_2$ rather than to cross the barrier~$f(c)-f(x_1)$ to  exit  through~$z_1$. Moreover, it can be proved that the quasi-stationary distribution~$\nu_h$ concentrates in any neighborhood of~$x_1$ in the limit $h\to 0$, which explains why the law of~$X_{\tau_\Omega}$  when~$X_0\sim \nu_h$  also concentrates on~$z_2$.  Concerning the two questions raised in the previous section, this example indicates that in the small temperature regime, there exist cases for which the process~\eqref{eq.langevin}, starting from the global minimum of~$f$ in~$\Omega$ or from~$\nu_h$, leaves~$\Omega$ through a point which is not a global minimum of~$f|_{\partial \Omega }$. 
\medskip

\noindent
This example also suggests the following. If one wants the law of~$X_{\tau_\Omega}$ to concentrate when~$h\to 0$ on points in $\pa \Omega$ where~$f$ attains its minimum, one should exclude cases when the largest timescales for the
diffusion process in~$\Omega$ are not related to energetic barriers involving points of~$\pa \Omega$  where~$f|_{\pa \Omega}$  attains its minimum. 
In  order to exclude such cases, we will assume in the following  that  the closure of each of  the connected components of~$\{f<\min_{\pa \Omega}f\}$  intersects~$\pa \Omega$.    \\

\noindent
Notice that if one modifies the function $f$ in the vicinity of $z_1$ such that $\pa_nf(z_1)>0$ and ${\rm argmin}_{\overline \Omega}f=\{x_1\}$, $z_1$ is then a generalized order one saddle point and the previous conclusions remain unchanged. 
\medskip

\noindent
\textit{Example 2}. Let us  construct a one-dimensional example for which the concentration of the law of~$X_{\tau_{\Omega}}$ on ${\rm argmin}_{\pa \Omega}f$ is not the same starting from the global minima of~$f$ in~$\Omega$ or from the quasi-stationary distribution~$\nu_h$. For this purpose, let us consider $
 z_1>0, \ z_2:=-z_1, \ z=0$ 
and $f\in C^{\infty}([z_1,z_2],\mathbb R) $ such that
 \begin{equation}\label{h1}
f \text{ is a Morse and even function}, \  \{x\in [z_1,z_2], f'(x)=0\}=\{x_1,z,x_2\},
 \end{equation}
where
 \begin{equation}\label{h2}
z_1<x_1<z<x_2<z_2, \, f(z_1)=f(z_2)>f(x_1)=f(x_2), \, f(z_1)< f(z).
 \end{equation}
A function $f$  satisfying~\eqref{h1} and~\eqref{h2} is represented in~Figure~\ref{fig:ex2}. One has the following result.
%
\begin{figure}
\begin{center}
\begin{tikzpicture}[scale=0.8]
  \draw [domain=1.37:3.56*pi, scale=0.5, samples=300, smooth] plot (\x,{2*cos(\x r) });
    \draw[->] (-0.5,-1.6) -- (6.1,-1.6);
  \draw[->] (-0.5,-1.6) -- (6.1,-1.6);
  \draw [dashed] (-0.2,0.22)--(6.5,0.22); 
   \draw (8,0.23) node[]{$\{f=\min_{\pa \Omega }f\}$};
  \draw (0.7,-1.83) node[]{$z_1$};
  \draw (0.7,-1.6) node[]{$\bullet$};
  \draw (5.5,-1.83) node[]{$z_2$};
  \draw (5.5,-1.6) node[]{$\bullet$};
  \draw (3.2,-1.83) node[]{$z$};
  \draw (3.2,-1.6) node[]{$\bullet$};
   \draw (1.58,-1.83) node[]{$x_1$};

  \draw (1.58,-1.6) node[]{$\bullet$};
    \draw (4.7,-1.83) node[]{$x_2$};
  \draw (4.7,-1.6) node[]{$\bullet$};
 
    \end{tikzpicture}    
 
    \caption{One-dimensional example where~\eqref{h1} and~\eqref{h2} are satisfied.}
    \label{fig:ex2}
    \end{center}
\end{figure}
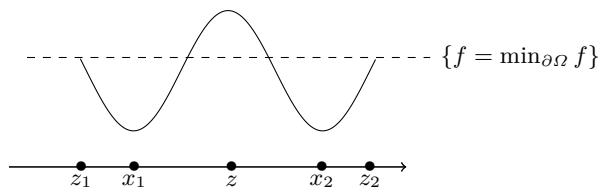
\begin{proposition}\label{pr.p-11}
Let $z_1>0$, $z_2:=-z_1$, $z=0$ and $f\in C^{\infty}([z_1,z_2],\mathbb R)$ which satisfies \eqref{h1} and~\eqref{h2}.   Then, one has for all $h>0$,
\begin{equation}\label{sortie-qsd}
\mathbb P_{\nu_h}  [X_{\tau_{(z_1,z_2)}}=z_1 ]=\frac 12   \ \text{ and }\  \mathbb P_{\nu_h}  [ X_{\tau_{(z_1,z_2)}}=z_2 ]=\frac 12  ,
\end{equation}
where  $\nu_h$ is the quasi-stationary distribution of the process~\eqref{eq.langevin} in $(z_1,z_2)$. Moreover, for all $x\in (z_1,z)$, there exists $c>0$ such that when $h\to 0$,
\begin{equation}\label{sortie-x1}
\mathbb P_{x}  [X_{\tau_{(z_1,z_2)}}=z_1 ]=1+O(e^{-\frac ch})\ \text{ and }\  \mathbb P_{x}  [ X_{\tau_{(z_1,z_2)}}=z_2 ]=O(e^{-\frac ch}),
\end{equation}
and for all $x\in (z,z_2)$, there exists $c>0$ such that when $h\to 0$
\begin{equation}\label{sortie-x2}
\mathbb P_{x}  [X_{\tau_{(z_1,z_2)}}=z_1 ]= O(e^{-\frac ch}) \ \text{ and }\  \mathbb P_{x}  [ X_{\tau_{(z_1,z_2)}}=z_2 ]=1+O(e^{-\frac ch}). 
\end{equation}
\end{proposition}
The asymptotic estimate~\eqref{sortie-qsd} is a consequence of the fact that $f$ is an even function (see~\cite[Section 1]{di-gesu-lelievre-le-peutrec-nectoux-16}). 
The asymptotic estimates~\eqref{sortie-x1} and~\eqref{sortie-x2} are proved exactly as Proposition~\ref{pr-ex1}, see~\cite[Section 1]{di-gesu-lelievre-le-peutrec-nectoux-16}. Let us also mention that Proposition~\ref{pr.p-11} is a consequence of the results~\cite{DoNe}. Concerning the two questions raised in the previous section,~Proposition~\ref{pr.p-11} shows that, when $f$ satisfies~\eqref{h1} and~\eqref{h2}, the concentration of the law of~$X_{\tau_{\Omega}}$ on $\{z_1,z_2\}$ is not the same starting from $x\in (z_1,z_2)\setminus\{z\}$ or from~$\nu_h$.  This is due to the fact that in this case the quasi-stationary distribution $\nu_h$ has an equal repartition in all disjoint neighboorhoods of $x_1$ and $x_2$, i.e.  for every  $(a_1,b_1)\subset (z_1,z)$ and  $(a_2,b_2)\subset (z,z_2)$ such that  $a_1<x_1<b_1$ and $a_2<x_2<b_2$,  it holds  for any $j\in \{1,2\}$, $\lim_{h\to 0}\int_{a_j}^{b_j}  \nu_h  = \frac{1}{  2 }$ (see~\cite{DoNe}). When $X_0=x\in (z_1,z_2)\setminus\{z\}$, the asymptotic estimates~\eqref{sortie-x1} and~\eqref{sortie-x2}  can be explained by the existence of a barrier $f(z)-f(x_1)$ which is larger than $f(z_1)-f(x_1)$.  
 In  order to exclude such cases, we will assume in the following  that there exists a connected component~$C$ of $\{f<\min_{\pa \Omega}f\}$, such that ${\rm argmin}_{\overline \Omega}f\subset C$.

\subsubsection{Main results on the exit point distribution.} 
\begin{figure}[!h]
   \begin{center}
  \begin{tikzpicture}[scale=0.7]
\tikzstyle{vertex}=[draw,circle,fill=black,minimum size=5pt,inner sep=0pt]

\draw (-7,2) ..controls  (-5.8,3.5).. (-5,2) ;
\draw (-5,2) ..controls  (-3.6,-0.9).. (-2,1 )  ;
\draw (-1.5,1 ) ..controls  (0,-0.9).. (1.8,2)  ;
\draw (-2,1 ) ..controls  (-1.75,1.33).. (-1.5,1)   ;
\draw (1.8,2) node[vertex,label=east: {$z_{1}$}](v){};

\draw (-7,2)  node[vertex,label=south: { }](v){}; 
 
\draw (-3.48,-0.3) node[vertex,label=south: {$x_{1}$}](v){};
\draw (-0.1,-0.3)  node[vertex,label=south: {$x_{2}$}](v){};

  \draw[dashed, <->]  (-4.85,2.05) -- (1.6,2.05);
 \draw (-1.6,2.5) node[]{\small{$C=\{f<\min_{\pa \Omega}f\}$}};

  \end{tikzpicture} 

 \caption{A one-dimensional example where   \textbf{[H-Morse]}  and \textbf{[H-Min]} are satisfied, the normal derivative of~$f$ changes sign  and the function $f$  has a saddle point in $\Omega$ higher than~$\min_{\pa \Omega}f$. In this example, $\{f<\min_{\pa \Omega}f\}$ is connected and thus $C=\{f<\min_{\pa \Omega}f\}$. Moreover, $\overline{ C }  \cap \pa \Omega=\{z_1\}$. }
 \label{fig:H0}
 \end{center}
 \end{figure}
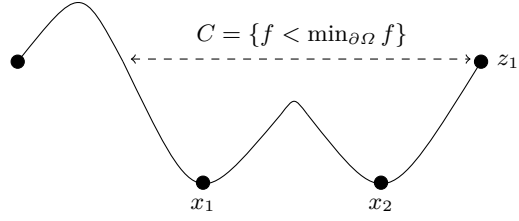

In this section, a simplified version of the results of~\cite{di-gesu-lelievre-le-peutrec-nectoux-16} is presented.  The aim is to exhibit  a simple geometric setting for which, on the one hand,   the law of $X_{\tau_\Omega}$ concentrates on the same points of $\pa \Omega$ when~$X_0\sim\nu_h$ or $X_0=x\in \Omega$ for some $x\in \{f<\min_{\pa \Omega}f\}$ and, on the other hand, this concentration occurs on generalized saddle points of~$f$ which belong to ${\rm argmin}_{\pa\Omega}f$. To this end, let us define the two following assumptions:
\begin{itemize}
\item \textbf{[H-Morse]} The function $f:\overline \Omega\to \mathbb R$  is $C^{\infty}$. The functions $f:\overline \Omega\to \mathbb R$  and~$f|_{\pa \Omega}$ are Morse  functions. Moreover, $\vert \nabla f\vert (x)\neq 0$ for all $x\in \pa \Omega$.\\
\item \textbf{[H-Min]}  The open set $\{f<\min_{\pa \Omega}f\}$ is nonempty, contains all the local minima of $f$ in $\Omega$  and the closure of each of the connected components of $\{f<\min_{\pa \Omega}f\}$  intersects $\pa \Omega$. Furthermore, there exists a connected component~$C$ of $\{f<\min_{\pa \Omega}f\}$ such that ${\rm argmin}_{\overline \Omega}f\subset C$. 
\end{itemize}
Notice that under \textbf{[H-Morse]} and \textbf{[H-Min]}, it holds $\min_{\pa \Omega}f >\min_{\overline  \Omega}f=\min_{ \Omega}f$. 
Under the assumptions \textbf{[H-Morse]} and \textbf{[H-Min]}, one defines the set of points $\{z_1,\ldots,z_{k_0}\}$ by
\begin{equation}\label{eq.PS-g}
\overline{ C }  \cap \pa \Omega=\{z_1,\ldots,z_{k_0}\}.
\end{equation}
 \begin{remark}
As already explained, the points $z_1,\ldots.,z_{k_0}$ are generalized saddle points of~$f$ on $\pa \Omega$ (see~\eqref{PS-gene}) since they satisfy
\begin{equation}\label{eq.inc}
\{z_1,\ldots,z_{k_0}\}\subset   \{z\in \pa \Omega, \, \pa_nf(z)>0\}\cap {\rm argmin}_{\pa \Omega} f.
\end{equation}
\end{remark}
\begin{remark}
Under \textbf{[H-Min]}, the normal derivative of~$f$ can change sign  and the function $f$  can have saddle points in $\Omega$ higher than~$\min_{\pa \Omega}f$, see for instance Figure~\ref{fig:H0}. 
\end{remark}

As shown in the following theorem,  the assumption \textbf{[H-Min]} ensures that the quasi-stationary distribution~$\nu_h$ concentrates in neighborhoods of the global minima of~$f$ in~$C$ and, starting from $x\in C$ or from $\nu_h$,  that the concentration of the law of $X_{\tau_\Omega}$  when $h\to 0$ occurs on the set of generalized saddle points  $\{z_1,\ldots,z_{k_0}\}$ (see~\eqref{eq.PS-g}). Notice that the assumption \textbf{[H-Min]} is not satisfied in the two examples given in the previous section (see Figures~\ref{fig:ex1-1} and~\ref{fig:ex2}). 

\begin{theorem}\label{thm2}
Let us assume that the hypotheses \textbf{[H-Morse]} and \textbf{[H-Min]}  are satisfied.  Let $\nu_h$  be the quasi-stationary distribution of the process~\eqref{eq.langevin} in $\Omega$ (see~\eqref{exp.qsd}).  Let $\mathcal V$ be an open subset of $\Omega$. 
 Then, if  $\mathcal V\cap {\rm argmin}_{C}f\neq \emptyset$,
 one has in the limit $h\to 0$:
$$
\nu_h\big(\mathcal  V\big)  = \frac{   \sum_{x\in \mathcal  V\cap {\rm argmin}_{C}f}  \big( {\rm det \ Hess } f   (x)   \big)^{-\frac12}  }{  \sum_{x\in {\rm argmin}_{C}f}  \big( {\rm det \ Hess } f   (x)   \big)^{-\frac12}  }\  \big(1+O(h) \big).
$$
When $\overline{\mathcal V}\cap  {\rm argmin}_{C}f= \emptyset$, 
 there exists $c>0$ such that  when $h \to 0$:
$$
\nu_h\big( \mathcal  V\big) =O\big ( e^{-\frac{c}{h}} \big).
$$
In addition, let  $F\in C^{\infty}(\pa \Omega,\mathbb R)$.  Then, when $h\to 0$: 
\begin{equation}\label{eq.FF}
\mathbb E_{\nu_h} \left [ F\left (X_{\tau_{\Omega}} \right )\right]=\sum_{i=1}^{k_0 } F(z_i)\,a_{i} +O(h^{\frac 14} ),
\end{equation}
where for $i\in \{1,\ldots,k_0\}$, 
\begin{equation}\label{eq.FFai}
 a_i=\frac{  \partial_nf(z_i)      }{  \sqrt{ {\rm det \ Hess } f\big|_{\partial \Omega}   (z_i) }  } \left (\sum \limits_{j=1}^{k_0} \frac{  \partial_nf(z_j)      }{  \sqrt{ {\rm det \ Hess } f\big|_{\partial \Omega}   (z_j) }  }\right)^{-1}. 
\end{equation}
Finally,~\eqref{eq.FF}    holds when  $X_0=x\in C $. 
\end{theorem}
\begin{remark}\label{re.genera}
In~\cite{di-gesu-lelievre-le-peutrec-nectoux-16}, one also gives sharp asymptotic estimates of~$\lambda_h$ and $\pa_n u_h$ in a more general setting than the one of Theorem~\ref{thm2} (for instance, we study the case when $f$ has local minima higher than $\min_{\pa \Omega}f$). However, in~\cite{di-gesu-lelievre-le-peutrec-nectoux-16}, we do not study the  precise asymptotic behaviour of $X_{\tau_{\Omega}}$ when $h\to 0$ near generalized saddle points~$z$  of $f$ on ${\pa \Omega}$ which are such that $f(z)>\min_{\pa \Omega}f$ as we did in~\cite{DLLN} (see Corollary~\ref{cc1}). Finally, in~\cite{di-gesu-lelievre-le-peutrec-nectoux-16}, the optimality of the remainder term $O(h^{\frac 14} )$ in~\eqref{eq.FF} is discussed and improved in some situations.  
\end{remark}
 \subsubsection{Ideas and sketch of the proof of Theorem~\ref{thm2}.}
 In this section, one gives the sketch of the proof of~\eqref{eq.FF} which is the main result of Theorem~\ref{thm2}. Recall that from~\eqref{expression-P},  for  $F\in C^{\infty}(\pa \Omega,\mathbb R)$
$$\mathbb E_{\nu_{h}}\big [ F(X_{\tau_{\Omega}}) \big]=-\frac {h}{2\lambda_h} \frac{\displaystyle \int_{\Sigma} F\,   \partial_n
  u_h  e^{-\frac{2}{h} f}   }{\displaystyle \int_\Omega u_h e^{-\frac{2}{h} f } },
$$
where $u_h$ is the eigenfunction associated with the principal eigenvalue $\lambda_h$ of~$L^{D,(0)}_{f,h}$. 
Therefore, to prove~\eqref{eq.FF}, one studies the asymptotic behaviour when $h\to 0$ of the following  quantities \begin{equation}\label{eq.qua}
\lambda_h,\ \pa_n u_h \text{ and  }  \int_\Omega u_h e^{-\frac{2}{h} f }.
\end{equation}
Under the assumptions \textbf{[H-Morse]} and \textbf{[H-Min]}, one defines 
$$m_0:=\text{Card } \Big( \{z\in \Omega,  z \text { is a local minimum of  }  f \}  \Big)$$
and  
\begin{align}
\nonumber
m_1:&= \text{Card } \Big ( \{ z \text { is a local minimum of  }  f|_{\pa \Omega}\}\cap\{z\in \pa \Omega, \, \pa_nf(z)>0\} \Big) \\
\label{eq.m1}
&\quad + \text{Card } \Big (  \{ z \text { is saddle point of }  f \}   \Big).
\end{align}
The integer $m_1$ is the number of generalized saddle points of $f$ in $\overline \Omega$ (see~\cite[Section 5.2]{helffer-nier-06}).
To study the asymptotic behaviour when $h\to 0$ of the quantities involved in~\eqref{eq.qua}, the starting  point is to again observe that $\nabla u_h$ is  solution to an eigenvalue problem   for the same eigenvalue~$\lambda_h$ (as already explained at the end of Section~\ref{sec.m1}). Indeed, $\nabla u_h$ is  solution to (see \eqref{eq:L1-intro})
\begin{equation}\label{eq:L1-intro-bis}
\left\{
\begin{aligned}
L^{(1)}_{f,h} \nabla u_h &= \lambda_h \nabla u_h \text{ in $\Omega$,}\\
\nabla_T u_h& = 0 \text{ on $\partial \Omega$,}\\
\left(-\frac{h}{2} {\rm div} + \nabla f \cdot \right) \nabla u_h & = 0 \text{ on $\partial \Omega$,}\\
\end{aligned}
\right.
\end{equation}
where we recall that $
L^{(1)}_{f,h}= -\frac{h}{2} \Delta + \nabla f \cdot \nabla + {\rm
  Hess} \, f$ is an operator acting on $1$-forms. Let us also recall that the operator  $L^{(1)}_{f,h}$ with tangential boundary conditions~\eqref{eq:L1-intro-bis} is denoted by~$L^{D,(1)}_{f,h}$. From~\eqref{eq:L1-intro-bis},  $\nabla u_h$ is an eigenform of $L^{D,(1)}_{f,h}$ associated with~$\lambda_h$. \\
The second ingredient is the following result: under the assumptions   \textbf{[H-Morse]} and \textbf{[H-Min]} and when $h\to 0$, the operator~$L^{D, (0)}_{f,h}$ has exactly~$m_0$ eigenvalues smaller than~$ \frac{\sqrt h}{2}$   and~$L^{D,(1)}_{f,h}$ has exactly~$m_1$ eigenvalues smaller than~$\frac{\sqrt h}{2}$ (see~\cite[Chapter 3]{helffer-nier-06}). 
Actually, all theses small eigenvalues are exponentially small when $h\to 0$, i.e. they are all   $O\big(e^{-\frac ch}\big)$ for some $c>0$ independent of~$h$. In particular~$\lambda_h$ is an exponentially small eigenvalue of~$L^{D, (0)}_{f,h}$ and of~$L^{D,(1)}_{f,h}$.
Let us denote by~$\pi_h^{(0)}$ (resp. $\pi_h^{(1)}$)  the orthogonal projector in $L_w^2(\Omega)$ onto the~$m_0$ (resp. $m_1$) smallest eigenvalues of  $L^{D, (0)}_{f,h}$  (resp. $L^{D,(1)}_{f,h}$). Then, according to the foregoing, one has   when $h\to 0$:
$$
  \dim \Ran\, \pi_h^{(0)} = m_0,  \ \ \dim \Ran\, \pi_h^{(1)} =m_1
$$
and 
$$
\nabla u_h\in \Ran\, \pi_h^{(1)}.
$$
Let us now explain how we prove Theorem~\ref{thm2}. To this end, let us introduce the set of local minima of $f$ in $\Omega$, 
$${\mathsf U}_0^\Omega:=\{x\in \Omega,\, x \text{ is a local minimum of } f\},$$
and the set of generalized saddle points of  $f$ in $\overline \Omega$,
\begin{align*}
{\mathsf U}_1^{\overline \Omega}&= \Big (  \{ z \text { is a local minimum of }  f|_{\pa \Omega}\}\cap\{z\in \pa \Omega, \, \pa_nf(z)>0\}  \Big) \\
&\quad \bigcup \   \{ z \text { is a saddle point of }  f \}.
\end{align*}
Let us recall that $m_0=\text{Card } \big ( {\mathsf U}_0^{\Omega} \big)$ and, from~\eqref{eq.m1},  that $m_1=\text{Card } \big ( {\mathsf U}_1^{\overline \Omega} \big)$. 
The first step to prove Theorem~\ref{thm2} consists in constructing  two maps $\widetilde {\mathbf j}$ and $  {\mathbf j}$. The goal of the map $\mathbf j$ is to associate each local minimum $x$ of~$f$ with a set of generalized saddle points $\mathbf j(x)\subset   {\mathsf U}_1^{\overline \Omega}$ such that 
$$\forall z,y\in \mathbf j(x), \ f(z)=f(y),$$
and such that, in the limit $h\to 0$, there exists at least one eigenvalue of $L^{D,(0)}_{f,h}$ whose exponential rate of decay is $2\big ( f(\mathbf j(x))-f(x)\big )$ i.e.
$$\exists \lambda \in \sigma\big (L^{D,(0)}_{f,h}\big) \, \text{ such that }\, \lim_{h\to 0}h\ln \lambda=-2 \big (f(\mathbf j(x))-f(x) \big ).$$
The aim of the map $\widetilde {\mathbf j}$ is to associate  each local minimum $x$ of~$f$ with the connected component of $\{f<f(\mathbf j(x))\}$ which contains $x$.   \medskip
 
 \noindent
The second step consists in constructing bases of $\Ran\, \pi_h^{(0)}$ and   $\Ran\, \pi_h^{(1)}$. To this end, one constructs two families of quasi-modes, denoted by $(\widetilde u_k)_{k\in \{1,\ldots,m_0\}}$ and $(\widetilde \psi_j)_{j\in \{1,\ldots,m_1\}}$,  which are then respectively  projected onto $\Ran\, \pi_h^{(0)}$ and  $\Ran\, \pi_h^{(1)}$.  To construct the family of  $1$-forms $(\widetilde \psi_j)_{j\in \{1,\ldots,m_1\}}$, we proceed as follows. For each saddle point~$z$ of $f$ in $\Omega$, following the procedure of~\cite{helffer-klein-nier-04}, one constructs a $1$-form supported in a neighboorhood of~$z$ in $\Omega$. For a local minimum~$z$ of $f|_{\partial \Omega }$  such that $\pa_n f(z)>0$, one constructs a $1$-form supported in a neighboorhood of~$z$ in $\overline \Omega$ as made in~\cite{helffer-nier-06}. To construct the family of functions $(\widetilde u_k)_{k\in \{1,\ldots,m_0\}}$,  one constructs for each local minimum $x$ of $f$ a smooth function whose support is almost $\widetilde {\mathbf j}(x)$ (this construction is close to the one made in~\cite{helffer-klein-nier-04,helffer-nier-06,HeHiSj,le-peutrec-10,michel2017small}).  
\medskip
 
 \noindent The next step consists in finding a sharp asymptotic equivalent for $\lambda_h$ when $h\to 0$. The quantity $\frac 2h  \lambda_h$ equals the square of the smallest singular values of the finite dimensional  operator  
$$\nabla: \Ran\, \pi_h^{(0)}\to \Ran\, \pi_h^{(1)}.$$
To study the asymptotic behaviour when $h\to 0$ of this smallest singular value, one uses the bases of $\Ran\, \pi_h^{(0)}$ and $\Ran\, \pi_h^{(1)}$ which have been constructed previously. 
The analysis of this finite dimensional problem is inspired by~\cite{HeHiSj} and also yields the asymptotic equivalent of $\int_{\Omega} u_h\,e^{-\frac 2h f}$ when $h\to 0$.\\
 
 \noindent
Then, we study the asymptotic behaviour of the normal derivative of~$u_h$ on~$\pa \Omega$ when $h\to 0$ to deduce  that the law of $X_{\tau_\Omega}$ concentrates when $h\to 0$ on  $\overline C \cap \pa \Omega=\{z_1,\ldots,z_{k_0}\}$ when $X_0\sim \nu_h$. \\
Lastly, one proves  ``leveling" results  on the function
$$x\mapsto \mathbb E_x[F(X_{\tau_\Omega})]$$
to obtain that when $X_0=x\in C$,  the law of~$X_{\tau_\Omega}$ also  concentrates when $h\to 0$  on $\{z_1,\ldots,z_{k_0}\}$.
\medskip
 
 \noindent
To conclude, the main results of~\cite{di-gesu-lelievre-le-peutrec-nectoux-16} are the following:
\begin{enumerate}
\item One uses techniques from semi-classical analysis   to study the asymptotic behaviours of  $\lambda_h$ and $\pa_n u_h$ when $h\to 0$,  and then, the concentration of the law of  $X_{\tau_\Omega}$  on  a subset of ${\rm argmin}_{\pa \Omega}f$ when $X_0\sim \nu_h$.

\item One identifies the points of  ${\rm argmin}_{\pa \Omega}f$  where the law of  $X_{\tau_\Omega}$  concentrates when $X_0\sim \nu_h$: this set of points is $\{z_1,\ldots,z_{k_0}\}$. Moreover, explicit formulas for their relative probabilities  are given (see indeed~\eqref{eq.FFai}) as well as precise remainder terms.

\item  One extends the previous results on the law of $X_{\tau_\Omega}$ to a deterministic initial condition in $\Omega$: $X_0=x$ where $x\in C$.

\item  \begin{sloppypar} These results hold under weak assumptions on the function $f$  and one-dimensional examples are given  to explain  why the geometric assumptions are needed to get them. \end{sloppypar}
\end{enumerate}
 \noindent
 \textbf{Conclusion.} We presented recent results which justify the use of a kinetic Monte Carlo model parametrized by Eyring-Kramers formulas to model the exit event from a metastable state $\Omega$ for the overdamped Langevin dynamics~\eqref{eq.langevin}. Our analysis is for the moment limited to  situations where $\vert \nabla f\vert \neq 0$ on $\pa \Omega$, which does not allow to consider order one saddle points on $\pa \Omega$. 
The extensions of~\cite{DLLN} and~\cite{di-gesu-lelievre-le-peutrec-nectoux-16} which are currently under study are the following:  the case when~$f$ has saddle points on $\pa \Omega$ and the case when the diffusion process  $X_t=(q_t,p_t)$ is solution to the Langevin stochastic differential equation
$$
\left\{
\begin{aligned}
 dq_t &= p_t dt,\\
dp_t&=-\nabla f(q_t)dt - \gamma \, p_t dt + \sqrt{h\gamma}\,  dB_t,
\end{aligned}
\right.
$$
where $(q_t,p_t)\in \Omega\times \mathbb R^d$, $\Omega$ being a bounded open subset of $\mathbb R^d$.



%
%
\medskip

\noindent
\begin{sloppypar}
\noindent
{\small\textbf{Acknowledgements.}  This work is supported by the European Research Council under the European
Union's Seventh Framework Programme (FP/2007-2013) / ERC Grant Agreement
number 614492. }
\end{sloppypar}
\bibliographystyle{plain}
 \bibliography{ihp_proc}

\begin{thebibliography}{10}

\bibitem{aristoff-lelievre-14}
D.~Aristoff and T.~Leli{\`e}vre.
\newblock Mathematical analysis of temperature accelerated dynamics.
\newblock {\em Multiscale Model. Simul.}, 12(1):290--317, 2014.

\bibitem{berglund-13}
N.~Berglund.
\newblock Kramers' law: Validity, derivations and generalisations.
\newblock {\em Markov Processes Relat. Fields}, 19:459--490, 2013.

\bibitem{berglund2014noise}
N.~Berglund.
\newblock Noise-induced phase slips, log-periodic oscillations, and the
  {G}umbel distribution.
\newblock {\em Markov Processes Relat. Fields}, 22:467--505, 2016.

\bibitem{berglund2004noise2}
N.~Berglund and B.~Gentz.
\newblock On the noise-induced passage through an unstable periodic orbit {I}:
  Two-level model.
\newblock {\em Journal of statistical physics}, 114(5-6):1577--1618, 2004.

\bibitem{berglund2014noise3}
N.~Berglund and B.~Gentz.
\newblock On the noise-induced passage through an unstable periodic orbit {II}:
  General case.
\newblock {\em SIAM Journal on Mathematical Analysis}, 46(1):310--352, 2014.

\bibitem{bovier-eckhoff-gayrard-klein-04}
A.~Bovier, M.~Eckhoff, V.~Gayrard, and M.~Klein.
\newblock Metastability in reversible diffusion processes. {I}. sharp
  asymptotics for capacities and exit times.
\newblock {\em J. Eur. Math. Soc. (JEMS)}, 6:399--424, 2004.

\bibitem{bovier-gayrard-klein-05}
A.~Bovier, V.~Gayrard, and M.~Klein.
\newblock Metastability in reversible diffusion processes. {II}. precise
  asymptotics for small eigenvalues.
\newblock {\em J. Eur. Math. Soc. (JEMS)}, 7:69--99, 2005.

\bibitem{cameron-14b}
M.~Cameron.
\newblock Computing the asymptotic spectrum for networks representing energy
  landscapes using the minimum spanning tree.
\newblock {\em Netw. Heterog. Media}, 9(3):383--416, 2014.

\bibitem{champagnat2015quasi}
N.~Champagnat and D.~Villemonais.
\newblock Quasi-stationary distribution for multi-dimensional birth and death
  processes conditioned to survival of all coordinates.
\newblock {\em arXiv preprint arXiv:1508.03161}, 2015.

\bibitem{chandrasekhar1943stochastic}
S.~Chandrasekhar.
\newblock Stochastic problems in physics and astronomy.
\newblock {\em Reviews of modern physics}, 15(1):1, 1943.

\bibitem{davies-82c}
E.B. Davies.
\newblock Dynamical stability of metastable states.
\newblock {\em J. Funct. Anal.}, 46(3):373--386, 1982.

\bibitem{davies-82a}
E.B. Davies.
\newblock Metastable states of symmetric {M}arkov semigroups {I}.
\newblock {\em Proc. London Math. Soc.}, 45(3):133--150, 1982.

\bibitem{davies-82b}
E.B. Davies.
\newblock Metastable states of symmetric {M}arkov semigroups {II}.
\newblock {\em J. London Math. Soc.}, 26(3):541--556, 1982.

\bibitem{day-83}
M.V. Day.
\newblock On the exponential exit law in the small parameter exit problem.
\newblock {\em Stochastics}, 8(4):297--323, 1983.

\bibitem{day1984a}
M.V. Day.
\newblock On the asymptotic relation between equilibrium density and exit
  measure in the exit problem.
\newblock {\em Stochastics: an international journal of probability and
  stochastic processes}, 12(3-4):303--330, 1984.

\bibitem{day1987r}
M.V. Day.
\newblock Recent progress on the small parameter exit problem.
\newblock {\em Stochastics: An International Journal of Probability and
  Stochastic Processes}, 20(2):121--150, 1987.

\bibitem{day1992}
M.V. Day.
\newblock Conditional exits for small noise diffusions with characteristic
  boundary.
\newblock {\em The Annals of Probability}, pages 1385--1419, 1992.

\bibitem{day1996}
M.V. Day.
\newblock Exit cycling for the van der {P}ol oscillator and quasipotential
  calculations.
\newblock {\em Journal of Dynamics and Differential Equations}, 8(4):573--601,
  1996.

\bibitem{day-99}
M.V. Day.
\newblock {\em Mathematical Approaches to the Problem of Noise-Induced Exit},
  pages 269--287.
\newblock Birkh{\"a}user, 1999.

\bibitem{devinatz-friedman-78a}
A.~Devinatz and A.~Friedman.
\newblock Asymptotic behavior of the principal eigenfunction for a singularly
  perturbed dirichlet problem.
\newblock {\em Indiana Univ. Math. J.}, 27:143--157, 1978.

\bibitem{devinatz-friedman-78b}
A.~Devinatz and A.~Friedman.
\newblock The asymptotic behavior of the solution of a singularly perturbed
  dirichlet problem.
\newblock {\em Indiana Univ. Math. J.}, 27(3):527--537, 1978.

\bibitem{DLLN}
G.~Di~Ges\`u, T.~Leli\`evre, D.~Le~Peutrec, and B.~Nectoux.
\newblock Sharp asymptotics of the first exit point density, 2017.
\newblock https://arxiv.org/pdf/1706.08728.pdf.

\bibitem{di-gesu-lelievre-le-peutrec-nectoux-16}
G.~Di~Ges\`u, T.~Leli\`evre, D.~Le~Peutrec, and B.~Nectoux.
\newblock The exit from a metastable state: concentration of the exit point
  distribution on the low energy saddle points, 2018.
\newblock In preparation.

\bibitem{eckhoff-05}
M.~Eckhoff.
\newblock Precise asymptotics of small eigenvalues of reversible diffusions in
  the metastable regime.
\newblock {\em Ann. Probab.}, 33(1):244--299, 2005.

\bibitem{fan-yip-yildiz-14}
Y.~Fan, S.~Yip, and B.~Yildiz.
\newblock Autonomous basin climbing method with sampling of multiple transition
  pathways: application to anisotropic diffusion of point defects in hcp {Zr}.
\newblock {\em Journal of Physics: Condensed Matter}, 26:365402, 2014.

\bibitem{freidlin-wentzell-84}
M.I. Freidlin and A.D. Wentzell.
\newblock {\em Random Perturbations of Dynamical Systems}.
\newblock Springer-Verlag, 1984.

\bibitem{galves-olivieri-vares-87}
A.~Galves, E.~Olivieri, and M.E. Vares.
\newblock Metastability for a class of dynamical systems subject to small
  random perturbations.
\newblock {\em The Annals of Probability}, 1288--1305, 1987.

\bibitem{goldshtein-mitrea-mitrea-11}
V.~Gol'dshtein, I.~Mitrea, and M.~Mitrea.
\newblock Hodge decompositions with mixed boundary conditions and applications
  to partial differential equations on {L}ipschitz manifolds.
\newblock {\em Journal of Mathematical Sciences}, 172(3):347--400, 2011.

\bibitem{hanggi-talkner-barkovec-90}
P.~H\"anggi, P.~Talkner, and M.~Borkovec.
\newblock Reaction-rate theory: fifty years after {K}ramers.
\newblock {\em Reviews of Modern Physics}, 62(2):251--342, 1990.

\bibitem{helffer-klein-nier-04}
B.~Helffer, M.~Klein, and F.~Nier.
\newblock Quantitative analysis of metastability in reversible diffusion
  processes via a {W}itten complex approach.
\newblock {\em Mat. Contemp.}, 26:41--85, 2004.

\bibitem{helffer-nier-06}
B.~Helffer and F.~Nier.
\newblock Quantitative analysis of metastability in reversible diffusion
  processes via a {W}itten complex approach: the case with boundary.
\newblock {\em M{\'e}moire de la Soci{\'e}t{\'e} math{\'e}matique de France},
  (105):1--89, 2006.

\bibitem{helffer-sjostrand-84}
B.~Helffer and J.~Sj\"ostrand.
\newblock Multiple wells in the semi-classical limit {I}.
\newblock {\em Comm. Partial Differential Equations}, 9(4):337--408, 1984.

\bibitem{helffer-sjostrand-85b}
B.~Helffer and J.~Sj{\"o}strand.
\newblock Multiple wells in the semi-classical limit {III}-{I}nteraction
  through non-resonant wells.
\newblock {\em Mathematische Nachrichten}, 124(1):263--313, 1985.

\bibitem{helffer-sjostrand-85a}
B.~Helffer and J.~Sj{\"o}strand.
\newblock Puits multiples en limite semi-classique. {II}. {I}nteraction
  mol{\'e}culaire. {S}ym{\'e}tries. {P}erturbation.
\newblock {\em Annales de l'IHP Physique th{\'e}orique}, 42(2):127--212, 1985.

\bibitem{HeSj4}
B.~Helffer and J.~Sj{\"o}strand.
\newblock Puits multiples en mecanique semi-classique iv etude du complexe de
  witten.
\newblock {\em Comm. Partial Differential Equations}, 10(3):245--340, 1985.

\bibitem{HeHiSj}
F.~H\'erau, M.~Hitrik, and J.~Sj\"ostrand.
\newblock Tunnel effect and symmetries for {K}ramers-{F}okker-{P}lanck type
  operators.
\newblock {\em J. Inst. Math. Jussieu}, 10(3):567--634, 2011.

\bibitem{holley-kusuoka-stroock-89}
R.A. Holley, S.~Kusuoka, and D.W. Stroock.
\newblock Asymptotics of the spectral gap with applications to the theory of
  simulated annealing.
\newblock {\em J. Funct. Anal.}, 83(2):333--347, 1989.

\bibitem{jakab-mitrea-mitrea-09}
T.~Jakab, I.~Mitrea, and M.~Mitrea.
\newblock On the regularity of differential forms satisfying mixed boundary
  conditions in a class of {L}ipschitz domains.
\newblock {\em Indiana Univ. Math. J.}, 58(5):2043--2071, 2009.

\bibitem{kamin-78}
S.~Kamin.
\newblock Elliptic perturbation of a first order operator with a singular point
  of attracting type.
\newblock {\em Indiana Univ. Math. J.}, 27(6):935--952, 1978.

\bibitem{kamin1979elliptic}
S.~Kamin.
\newblock On elliptic singular perturbation problems with turning points.
\newblock {\em SIAM Journal on Mathematical Analysis}, 10(3):447--455, 1979.

\bibitem{kipnis-newman-85}
C.~Kipnis and C.M. Newman.
\newblock The metastable behavior of infrequently observed, weakly random,
  one-dimensional diffusion processes.
\newblock {\em SIAM J. Appl. Math.}, 45(6):972--982, 1985.

\bibitem{kramers-40}
H.A. Kramers.
\newblock Brownian motion in a field of force and the diffusion model of
  chemical reactions.
\newblock {\em Physica}, 7(4):284--304, 1940.

\bibitem{Kra}
H.A. Kramers.
\newblock Brownian motion in a field of force and the diffusion model of
  chemical reactions.
\newblock {\em Physica}, 7(4):284--304, 1940.

\bibitem{landim2017dirichlet}
C.~Landim, M.~Mariani, and I.~Seo.
\newblock A {D}irichlet and a {T}homson principle for non-selfadjoint elliptic
  operators, metastability in non-reversible diffusion processes.
\newblock {\em arXiv preprint arXiv:1701.00985}, 2017.

\bibitem{le-bris-lelievre-luskin-perez-12}
C.~Le~Bris, T.~Leli{\`e}vre, M.~Luskin, and D.~Perez.
\newblock A mathematical formalization of the parallel replica dynamics.
\newblock {\em Monte Carlo Methods and Applications}, 18(2):119--146, 2012.

\bibitem{le-peutrec-10}
D.~Le~Peutrec.
\newblock Small eigenvalues of the {N}eumann realization of the semiclassical
  {W}itten {L}aplacian.
\newblock {\em Ann. Fac. Sci. Toulouse Math. (6)}, 19(3-4):735--809, 2010.

\bibitem{DoNe}
D.~Le~Peutrec and B.~Nectoux.
\newblock Repartition of the quasi-stationary distribution and first exit point
  density for a double-well potential, 2018.
\newblock In preparation.

\bibitem{marcelin-15}
R.~Marcelin.
\newblock Contribution \`a l'\'etude de la cin\'etique physico-chimique.
\newblock {\em Ann. Physique}, 3:120--231, 1915.

\bibitem{mathieu-94}
P.~Mathieu.
\newblock Zero white noise limit through {D}irichlet forms, with application to
  diffusions in a random medium.
\newblock {\em Probability Theory and Related Fields}, 99(4):549--580, 1994.

\bibitem{mathieu-95}
P.~Mathieu.
\newblock Spectra, exit times and long time asymptotics in the zero-white-noise
  limit.
\newblock {\em Stochastics}, 55(1-2):1--20, 1995.

\bibitem{MS77}
B.J. Matkowsky and Z.~Schuss.
\newblock The exit problem for randomly perturbed dynamical systems.
\newblock {\em SIAM J. Appl. Math.}, 33(2):365--382, 1977.

\bibitem{MaSc}
B.J. Matkowsky and Z.~Schuss.
\newblock The exit problem: a new approach to diffusion across potential
  barriers.
\newblock {\em SIAM J. Appl. Math.}, 36(3):604--623, 1979.

\bibitem{MATK1982}
B.J. Matkowsky, Z.~Schuss, and E.~Ben-Jacob.
\newblock A singular perturbation approach to kramers diffusion problem.
\newblock {\em SIAM J. Appl. Math.}, 42(4):835--849, 1982.

\bibitem{michel2017small}
L.~Michel.
\newblock About small eigenvalues of {W}itten laplacian.
\newblock {\em arXiv preprint arXiv:1702.01837}, 2017.

\bibitem{miclo-95}
L.~Miclo.
\newblock Comportement de spectres d'op\'erateurs de {S}chr\"odinger \`a basse
  temp\'erature.
\newblock {\em Bulletin des sciences math{\'e}matiques}, 119(6):529--554, 1995.

\bibitem{schuss90}
T.~Naeh, M.M. Klosek, B.J. Matkowsky, and Z.~Schuss.
\newblock A direct approach to the exit problem.
\newblock {\em SIAM J. Appl. Math.}, 50(2):595--627, 1990.

\bibitem{BN2017}
B.~Nectoux.
\newblock {\em Analyse spectrale et analyse semi-classique pour la
  m{\'e}tastabilit{\'e} en dynamique mol{\'e}culaire}.
\newblock PhD thesis, Universit{\'e} {P}aris {E}st, 2017.

\bibitem{perthame-90}
B.~Perthame.
\newblock Perturbed dynamical systems with an attracting singularity and weak
  viscosity limits in hamilton-jacobi equations.
\newblock {\em Transactions of the American Mathematical Society},
  317(2):723--748, 1990.

\bibitem{schilder1966some}
M.~Schilder.
\newblock Some asymptotic formulas for {W}iener integrals.
\newblock {\em Transactions of the American Mathematical Society},
  125(1):63--85, 1966.

\bibitem{schuette-98}
C.~Sch\"utte.
\newblock Conformational dynamics: modelling, theory, algorithm and application
  to biomolecules, 1998.
\newblock Habilitation dissertation, Free University Berlin.

\bibitem{schuette-sarich-13}
C.~Sch\"utte and M.~Sarich.
\newblock {\em Metastability and {M}arkov state models in molecular dynamics},
  volume~24 of {\em Courant Lecture Notes}.
\newblock American Mathematical Society, 2013.

\bibitem{sorensen-voter-00}
M.R. Sorensen and A.F. Voter.
\newblock Temperature-accelerated dynamics for simulation of infrequent events.
\newblock {\em J. Chem. Phys.}, 112(21):9599--9606, 2000.

\bibitem{sugiura-95}
M.~Sugiura.
\newblock Metastable behaviors of diffusion processes with small parameter.
\newblock {\em Journal of the Mathematical Society of Japan}, 47(4):755--788,
  1995.

\bibitem{vineyard-57}
G.H. Vineyard.
\newblock Frequency factors and isotope effects in solid state rate processes.
\newblock {\em Journal of Physics and Chemistry of Solids}, 3(1):121--127,
  1957.

\bibitem{voter-97}
A.F. Voter.
\newblock A method for accelerating the molecular dynamics simulation of
  infrequent events.
\newblock {\em J. Chem. Phys.}, 106(11):4665--4677, 1997.

\bibitem{voter-98}
A.F. Voter.
\newblock Parallel replica method for dynamics of infrequent events.
\newblock {\em Phys. Rev. B}, 57(22):R13 985, 1998.

\bibitem{voter-05}
A.F. Voter.
\newblock {\em Radiation Effects in Solids}, chapter Introduction to the
  Kinetic Monte Carlo Method.
\newblock Springer, NATO Publishing Unit, 2005.

\bibitem{wales-03}
D.J. Wales.
\newblock {\em Energy landscapes}.
\newblock Cambridge University Press, 2003.

\end{thebibliography}

%
%
%
%
%
%
%
%
\end{document}